\documentclass[10pt]{amsart}
\usepackage{graphicx,color}
\usepackage{color}
\usepackage{url}
\newtheorem{theorem}{Theorem}[section]
\newtheorem{lemma}[theorem]{Lemma}           
\newtheorem{cor}[theorem]{Corollary}

\theoremstyle{definition}
\newtheorem{definition}[theorem]{Definition}

\theoremstyle{remark}

\numberwithin{equation}{section}

\title[Discrete Schr{\"o}dinger Operators]{Inverse scattering at a fixed energy for Discrete Schr{\"o}dinger Operators on the square lattice}
\author{Hiroshi ISOZAKI}
\address{Institute of Mathematics,
University of Tsukuba
\\
Tsukuba, 305-8571, JAPAN 
\\ isozakih@math.tsukuba.ac.jp, hmorioka@math.tsukuba.ac.jp}
\author{Hisashi Morioka}

\begin{document}
\maketitle

\begin{abstract}
We study an inverse scattering problem for the discrete Schr{\"o}dinger operator  on the square lattice ${\bf Z}^d$, $d \geq 2$, with compactly supported potential. We show that the potential is uniquely reconstructed from a scattering matrix for a fixed energy.
\end{abstract}

%%%%%%%% Section 1 %%%%%%%%%%%%%%%%%%%

\section{Introduction}

%%%%%%%%%%%%%%%% subsection 1.1 %%%%%%%%%%%%%%%%%

\subsection{Inverse scattering}
Let ${\bf Z}^d = \{n = (n_1,\cdots,n_d) \, ; \, n_i \in {\bf Z}, \ 1 \leq i \leq d\}$ be the square lattice, and $e_{1} = (1,0,\cdots,0), \cdots, e_{d} = (0,\cdots,0,1)$ the standard
basis of ${\bf Z}^d$. Throughout the paper, we shall assume that $d \geq 2$. The  Schr{\"o}dinger operator
${\widehat H}$ on ${\bf Z}^d$ is defined by
\begin{equation}
{\widehat H} = {\widehat H}_{0} + {\widehat V},
\nonumber
\end{equation}
where for ${\widehat f} = \{\widehat f(n)\}_{n\in{\bf Z}^d} \in \ell^{2}({\bf Z}^d)$ and $n \in {\bf Z}^d$
\begin{equation}
\big({\widehat H}_{0}\widehat f\big)(n) = - \frac{1}{4}\sum_{j=1}^{d} \big\{ {\widehat f}(n + e_{j}) +
{\widehat f}(n - e_{j}) \big\} + \frac{d}{2}\widehat f(n),
\nonumber
\end{equation}
\begin{equation}
({\widehat V}{\widehat f})(n) = \widehat V(n){\widehat f}(n).
\nonumber
\end{equation}
We impose the following assumption on $\widehat V$:

\medskip
\noindent
{\bf (A)}  $\ \widehat V$ {\it is real-valued, and 
$\widehat V(n) = 0$ except for a finite number of} $n$.

\medskip
Under this assumption, $\sigma(\widehat H_0) = \sigma_{ess}(\widehat H) = [0,d]$, and  the wave operators
\begin{equation}
\widehat W^{(\pm)} = {\mathop{\rm s-lim}_{t\to\pm \infty}}\,e^{it\widehat H}e^{-it\widehat H_0} \quad ( {\rm in} \ \ell^2({\bf Z}^d))
\label{S1WaveOp}
\end{equation}
exist and are asymptotically complete, i.e. their ranges coincide with ${\mathcal H}_{ac}(\widehat H)$, the absolutely continuous subspace for $\widehat H$. Hence the scattering operator
\begin{equation}
\widehat S = \big(\widehat W^{(+)}\big)^{\ast}\widehat W^{(-)}
\end{equation}
is unitary. Associated with $\widehat H_0$, we have a unitary spectral representation 
$$
\widehat{\mathcal F}_0 : \ell^2({\bf Z}^d) \to L^2((0,d);L^2(M_{\lambda});d\lambda),
$$
 where
\begin{equation}
M_{\lambda} = \Big\{x \in {\bf T}^d \, ; \, d - \sum_{j=1}^d\cos x_j = 2\lambda \Big\},
\label{S1Mlambda}
\end{equation}
\begin{equation}
{\bf T}^d = {\bf R}^d/(2\pi {\bf Z})^d = [-\pi,\pi]^d.
\label{S1Torus}
\end{equation}
Then 
$\widehat{\mathcal F}_0\widehat S(\widehat{\mathcal F}_0)^{\ast}$ has the following direct integral representation
\begin{equation}
\widehat{\mathcal F}_0\widehat S(\widehat{\mathcal F}_0)^{\ast} = \int_{0}^d\oplus \, \mathcal S(\lambda)\, d\lambda.
\label{S1DirectInt}
\end{equation}
Here $\mathcal S(\lambda)$ is a unitary operator on $L^2(M_{\lambda})$, and is called the {\it S-matrix}.

Our main concern in this paper is the inverse scattering, i.e. reconstruction of the potential $\widehat V$ from the knowledge of the S-matrix. In \cite{IsKo} (see also \cite{Es}), it has been proven that given $\mathcal{S} (\lambda)$ for all energy $\lambda \in (0,d)\setminus{\bf Z}$, one can uniquely reconstruct the potential.

It is worthwhile to recall the case of the continuous model, i.e. the Schr{\"o}dinger operator $- \Delta + V(x)$ in $L^2({\bf R}^d)$. In this case, it is known that only one arbitrarily fixed energy $\lambda > 0$ is sufficient to reconstruct the compactly supported (and also exponentially decaying) potential $V(x)$ from the S-matrix $\mathcal S(\lambda)$. This was proved  for $d \geq 3$ in 1980's by  Sylvester-Uhlmann \cite{SyUh87}, Nachman \cite{Na1}, Khenkin-Novikov \cite{KheNo}. There are two methods. One way is applicable to the compactly supported potential and based on the equivalence of the S-matrix and the Dirichlet-Neumann map (called D-N map hereafter) for the boundary value problem in a bounded domain. The other way relies on Faddeev's theory for the multi-dimensional inverse scattering, in particular, on Feddeev's scattering amplitude, and allows exponentially decaying potentials. In both cases, Sylvester-Uhlmann's complex geometrical optics solutions to the Schr{\"o}dinger equation, or Faddeev's exponentially growing Green function played a crucial role. (See e.g. an expositiory article \cite{Is03}.) However, since both of these methods use the complex Born approximation, the case $d=2$ remained open rather long time. 
 Note that for the potential of the form coming from ellectric conductivities, the 2-dim. inverse scattering problem for a fixed energy was solved by Nachman \cite{Nach96}. See also \cite{Isa}.
Recently Bukhgeim \cite{Buk}  proved that, based on  Carleman estimates, the D-N map determines the potential for the 2-dim. boundary value problem. For the partial data problem, see \cite{ImUhYa}. This result can be applied to the inverse scattering and to derive an affirmative answer to the uniqueness of the potential  for given potential of fixed energy.

\subsection{Main result}
To study the inverse scattering from a fixed energy for the discrete model, we adopt the above-mentioned former approach. Namely, we assume that the potential is compactly supported, and derive the equivalence of the S-matrix and the D-N map in a bounded domain. 

We need to restrcit the energy in some interval. Let 
\begin{equation}
I_d = \left\{
\begin{array}{lc}
(0,1)\cup(0,1), & {\rm for} \quad d = 2, \\
(0,1/2)\cup(d - 1/2 ,d), &{\rm for} \quad d \geq 3.
\end{array}
\right.
\label{S1IntervalI}
\end{equation}
The following theorem is our main aim.

%%%%%%%%%%%%%%%%%% Theorem 1.1 %%%%%%%%%%%%%%%

\begin{theorem} \label{Main theorem}
Fix $\lambda\in I_d$ arbitrarily. Then from the S-matrix $\mathcal S(\lambda)$, 
one can uniquely reconstruct the potential $\widehat V$.
\end{theorem}

Our proof not only  states  the uniqueness, but also explains the procedure of the reconstuction of the potential.

\subsection{The plan of the proof}
After the preparation of basic spectral results in \S 2 and \S3, the first task is to relate the S-matrix with the far-field pattern  at infinity of the generalized eigenfunction of $\widehat H$. This is done in \S4 by observing the asymptotic expansion at infinity of the Green operator of $\widehat H$. 
In \S5, we introduce the radiation condition for the Helmholtz equation and prove the uniqueness theorem for the solution.
We then study the spectral theory for the exterior problem in \S 6, with the aid of which we obtain in \S 7 the equivalence of the S-matrix and the D-N map for a boundary value problem in a bounded domain. The potential is then reconstructed from the D-N map in \S 8 via a constructive procedure. 

Although the main stream of the proof is the same as the continuous case, we need to be careful about the difference in the case of the discrete model.
 The first one is the asymptotic expansion of the resolvent at infinity. This is based on the stationary phase method on the surface $M_{\lambda}$ defined by 
 (\ref{S1Mlambda}), which is not strictly convex in general. This is  the reason we restrict the energy on $I_d$. 
The second one, which is more serious, occurs when we compare the far-field patterns of solutions to Schr{\"o}dinger equations in the whole space with those of the exterior domain. We need a Rellich type theorem (see Theorem \ref{Rellictype}) and a unique continuation property for the discrete Helmholtz equation, which do not seem to be well-known. 
However, the former's precursor has been given by Shaban-Vainberg \cite{Sha},  and the latter follows rather easily from it. As a byproduct, it proves the non-existence of embedded eigenvalues for $\widehat H$ (\cite{IsMo}).
We then go into the final step of computing the potential from the D-N map. In the continuous case, this is an elliptic Cauchy problem from the boundary, hence is ill-posed. However, in the discrete case, this is a finite dimensional problem, therefore a finite computational procedure. The whole proof does not depend on the space dimension. In contrust, it took a long time to get the 2-dim. result in the continuous case.

\subsection{Remarks for references}
There are important precursors of this paper. The work of Eskina \cite{Es}
 have already announced the result of the inverse scattering for discrete Schr{\"o}dinger operators. In particular, this paper stresses the effectiveness of several complex variables in the study of discrete Schr{\"o}dinger operators. Shaban-Vainberg \cite{Sha} studied the spectral theory of discrete Schr{\"o}dinger operators. They introduced the radiation condition, proved the limiting absorption principle, and derived the asymptotic expansion of the resolvent at infinity including the case of non-convex surface.

The computation of the D-N map for the discrete interior boundary value problem was done in the work of Oberlin \cite{Ob}. 
See also Curtis-Morrow \cite{Cu1} and Curtis-Mooers-Morrow \cite{Cu2}.

\subsection{Notation}
$C$'s denote various constants.
For any $x,y \in {\bf R}^d $, $x\cdot y= x_1 y_1 +\cdots +x_d y_d $ denotes the ordinary scalar product in the Euclidean space where $x_j $ and $y_j $ are $j$-th component of $x$ and $y$ respectively.
For any $x\in {\bf R}^d $, $|x| =( x \cdot x  )^{1/2} $ is the Euclidean norm.
Note that even for $n= (n_1 , \cdots , n_d ) \in {\bf Z}^d $, we use $ |n|=\big( \sum_{j=1}^d |n_j|^2 \big) ^{1/2} $.
For two Banach spaces $X$ and $Y$, $ \mathbf{B} (X,Y)$ denotes the totality of bounded operators from $X$ to $Y$.
For a self-adjoint operator $A$ on a Hilbert space, $\sigma (A)$, $\sigma _{ess} (A) $, $\sigma_{disc} (A)$, $\sigma_{ac} (A)$ and $\sigma_p (A)$ denote  its spectrum, essential spectrum, discrete spectrum, absolutely continuous spectrum and point spectrum, respectively. 
For a set $S$, $\, ^{\#} S $ denotes the number of elements in $S$.
We use the notation
$$
\langle t\rangle = (1 + t^2)^{1/2}, \quad t \in {\bf R} .
$$

\subsection{Acknowledgement}
The authors are indebted to Evgeny Korotyaev for useful discussions and  encouragements. The second author is supported by the Japan Society for the Promotion of Science under the Grant-in-Aid for Research Fellow (DC2) No. 23110.

%%%%%%%%%%%% %%%%%%%%%%%%% Section 2 %%%%%%%%%%%%%%

\section{Momentum representation}

%%%%%%%%%%%%%%%%%%%%%% subsection 2.1 %%%%%%%%%%%

\subsection{Discrete Fourier transform}
From the view point of dynamics on the lattice, the torus ${\bf T}^d$  in (\ref{S1Torus}) plays the role of momentum space. Let ${\mathcal U}$ be the unitary operator from
$\ell^{2}({\bf Z}^d)$ to $L^{2}({\bf T}^d)$ defined by
\begin{equation}
({\mathcal U}\,{\widehat f})(x) = (2\pi)^{-d/2}\sum_{n\in{\bf Z}^d}{\widehat f}(n)
e^{-in\cdot x}.
\nonumber
\end{equation}
Using this discrete Fourier transformation, the Hamiltonian $\widehat H$ is represented by
\begin{equation}
H= {\mathcal U}\,{\widehat H}\,{\mathcal U}^{\ast} = H_0 + V, \quad
H_0 = {\mathcal U}\, {\widehat H_0}\, {\mathcal U}^{\ast}, \quad
V = {\mathcal U}\, {\widehat V}\, {\mathcal U}^{\ast},
\nonumber
\end{equation}
where $H_0$ is the multiplcation operator:
\begin{equation}
H_0 = \frac{1}{2}\Big(d - \sum_{j=1}^d\cos x_j\Big) = :h(x),
\label{S2H0cosine}
\end{equation}
and $V$ is the convolution operator
\begin{equation}
(Vu)(x) = (2\pi)^{-d/2}\int_{{\bf T}^d}V(x-y)u(y)dy, \quad
V(x) = (2\pi)^{-d/2}\sum_{n\in{\bf Z}^d}\widehat V(n)e^{-in\cdot x}.
\nonumber
\end{equation}

%%%%%%%%%%%%%%%%%%%%%%%% subsection 2.2 %%%%%%%%%%%%%%%%%%%

\subsection{Sobolev and Besov spaces}
We define operators $\widehat N_j$ and $N_j$ by
$$
\big(\widehat N_j\widehat f)(n) = n_j\widehat f(n),
\quad
N_j = \mathcal U\widehat N_j\mathcal U^{\ast} = i\frac{\partial}{\partial x_j}.
$$
We put $N = (N_1,\cdots,N_d)$, and let $N^2$ be the self-adjont operator defined by
\begin{equation}
N^2 = \sum_{j=1}^dN_j^2 =  - \Delta, \quad {\rm on} \quad {\bf T}^d,
\nonumber
\end{equation}
where $\Delta$ denotes the Laplacian on ${\bf T}^d = [-\pi,\pi]^d$ with periodic boundary condition.  We put
\begin{equation}
|N| = \sqrt{N^2} = \sqrt{-\Delta}.
\nonumber
\end{equation}
For $s \in {\bf R}$,
let ${\mathcal H}^s$ be the completion of $D(|N|^s)$ with respect to the norm
$\|u\|_{s} = \|\langle N\rangle^{s}u\|$ :
\begin{equation}
{\mathcal H}^s = \{u \in {\mathcal D}^{\prime}({\bf T}^d)\, ; \, \|u\|_{s} = \|\langle N\rangle^{s}u\| <
\infty\},
\nonumber
\end{equation}
where $ \mathcal{D}' ({\bf T}^d ) $ denotes the space of distribution on ${\bf T}^d $.
Put $\mathcal H = \mathcal H^0 = L^2({\bf T}^d)$.

For a self-adjoint operator $T$, let $\chi(a \leq T < b)$ denote the operator $\chi_{I}(T)$, where $\chi_I(\lambda)$ is the characteristic function of the interval $I = [a,b)$. The operators $\chi(T < a)$ and $\chi(T \geq b)$ are defined similarly.  Using the series $\{r_j\}_{j=0}^{\infty}$ with $r_{-1} = 0$, $r_j = 2^j \ (j \geq 0)$, we define the Besov space $\mathcal B$ by
\begin{equation}
\mathcal B = \Big\{f \in {\mathcal H}\, ; \|f\|_{\mathcal B} = \sum_{j=0}^{\infty}r_j^{1/2}\|\chi(r_{j-1} \leq |N| < r_j)f\| < \infty\Big\}.
\nonumber
\end{equation}
Its dual space $\mathcal B^{\ast}$ is the completion of $\mathcal H$ by the following norm
\begin{equation}
\|u\|_{\mathcal B^{\ast}}= \sup_{j\geq 0}r_j^{-1/2}\|\chi(r_{j-1} \leq |N| < r_j)u\|.\nonumber
\end{equation}
The following Lemma \ref{S2sobolev_besov} is proved in the same way as in  \cite{AgHo76}.

%%%%%%%%%%%%%%%%%%%%%%%%%% Lemma 2.2 %%%%%%%%%%%%%%%%%%%%%%%%%

\begin{lemma}
(1) There exists a constant $C > 0$ such that
\begin{equation}
C^{-1}\|u\|_{\mathcal B^{\ast}}\leq \left(\sup_{R>1}\frac{1}{R}\|\chi(|N| < R)u\|^2\right)^{1/2}\leq C\|u\|_{\mathcal B^{\ast}}.
\nonumber
\end{equation}
Therefore, in the following, we use
\begin{equation}
\|u\|_{\mathcal B^{\ast}} = \left(\sup_{R>1}\frac{1}{R}\|\chi(|N| < R)u\|^2\right)^{1/2}
\nonumber
\end{equation}
as a norm on $\mathcal B^{\ast}$.\\
\noindent
(2) For $s > 1/2$, the following inclusion relations hold :
\begin{equation}
\mathcal H^{s} \subset \mathcal B \subset \mathcal H^{1/2} \subset \mathcal H \subset \mathcal H^{-1/2} \subset \mathcal B^{\ast} \subset \mathcal H^{-s}.
\nonumber
\end{equation}
\label{S2sobolev_besov}
\end{lemma}

\medskip
We also put $\widehat{\mathcal H} = \ell^2({\bf Z}^d)$, and define $\widehat{\mathcal H}^s$, $\widehat{\mathcal B}$, $\widehat{\mathcal B}^{\ast}$ by replacing $N$ by $\widehat N$. Note that $\widehat{\mathcal H}^s = \mathcal U^{\ast}\mathcal H^s$ and so on. In particular, 
Parseval's formula implies that
\begin{equation}
\|u\|_{\mathcal H^s}^2 = \|\widehat u\|_{\widehat{\mathcal H}^s}^2 = \sum_{n\in {\bf Z}^d}(1 + |n|^2)^s|\widehat u(n)|^2 ,
\nonumber
\end{equation}
\begin{equation}
\|u\|_{{\mathcal B}^{\ast}}^2 =  
\|\widehat u\|_{\widehat{\mathcal B}^{\ast}}^2
= 
\sup_{R>1}\frac{1}{R}\sum_{|n|<R}|\widehat u(n)|^2,
\nonumber
\end{equation}
$\widehat u(n)$ being the Fourier coefficient of $u(x)$.

%%%%%%%%%%%%%%%%%%%%%%%%%%%%% subsection 2.3 %%%%%%%%%%%%%%%%%%%%%

\subsection{Resolvent estimate}

%%%%%%%%%%%%%%%%%%%%% Lemma 2.1 %%%%%%%%%%%%%%%

\begin{lemma} 
(1) $\  \sigma(\widehat H_{0}) = \sigma_{ac} (\widehat H_{0}) = [0,d]$. \\
\noindent
(2) $\  $ $\sigma_{ess}(\widehat H) = [0,d], \quad
\sigma_{disc}(\widehat H) \subset {\bf R}\setminus[0,d]$. \\
\noindent
(3) $\ $ $ \sigma_{p}(\widehat H)\cap \big((0,d)\setminus {\bf Z}\big) = \emptyset$.
\end{lemma}

Proof. The assertions (1), (2) follow from (\ref{S2H0cosine}) and Weyl's theorem. The assertion (3) is proven in \cite{IsMo}. \qed

\medskip
Let $\widehat R(z) = (\widehat H - z)^{-1}$.

%%%%%%%%%%%%%%%%%%%% Theorem 2.3 %%%%%%%%%%%%%%%

\begin{theorem}

\noindent
(1) Let $s > 1/2$ and $\lambda \in (0,d)\setminus{\bf Z}$.  Then there exists a norm limit $\widehat  R(\lambda \pm i0) := \lim_{\epsilon\to0} \widehat R(\lambda \pm i\epsilon) \in {\bf B} ( \widehat{\mathcal{H}}^s ; \widehat{\mathcal{H}}^{-s} ) $. Moreover, we have
\begin{equation}
\sup_{\lambda \in J} \|\widehat R(\lambda \pm i0)\|_{{\bf B}(\widehat{\mathcal B};\widehat{\mathcal B}^{\ast})} < \infty.
\label{S2LAP}
\end{equation}
for any  compact interval $J$ in $(0,d)\setminus{\bf Z}$. The mapping $(0,d)\setminus{\bf Z} \ni \lambda \mapsto \widehat R(\lambda \pm i0)$ is norm continuous in ${\bf B}(\widehat{\mathcal H}^s;\widehat {\mathcal H}^{-s})$ and weakly continuous in ${\bf B}(\widehat {\mathcal B}\,; \widehat{\mathcal B}^{\ast})$.\\
\noindent
(2) $\widehat H$ has no singular continuous spectrum. \label{LAP}
\end{theorem}

For the proof of Theorem \ref{LAP}, see Lemma 2.5 and Theorem 2.6 of \cite{IsKo}. Note that 
\begin{equation}
\nabla h(x) = 0 \Longleftrightarrow h(x) \in \{0, 1, \cdots, d\}.
\label{S2nablah(x)=0}
\end{equation}
This is the reason why the set of thresholds $\{0, 1, \cdots, d\}$ appears.

%%%%%%%%%%%%%%%% subsection 3.2 %%%%%%%%%%%%%%%%%%

%%%%%%%%%%%%%%% Section 3 %%%%%%%%%%%%%%%%%

\section{Spectral representations and S-matrices}

%%%%%%%%%%%%%%%%% subsection 3.1 %%%%%%%%%%%%%%%
We recall spectral representations and S-matrices derived in \S 3 of \cite{IsKo}.
\subsection{Spectral representation on the torus}
We begin with the spectral representation in the momentum space.
Let us note
\begin{equation}
h(x) = \frac{1}{2}\Big(d - \sum_{j=1}^d\cos x_j\Big) = \sum_{j=1}^d\sin^2\left(\frac{x_j}{2}\right),
\nonumber
\end{equation}
which suggests that the variables $y = (y_1,\cdots,y_d) \in [-1,1]^d$:
\begin{equation}
y_j = \sin\frac{x_j}{2}, \quad x_j = 2\arcsin y_j
\nonumber
\end{equation}
are convenient to describe $H_0$. 
Note that for $\lambda \in ( 0 , d ) \setminus {\bf Z}$
\begin{equation}
x(\sqrt{\lambda}\theta) = \left(2\arcsin(\sqrt\lambda\theta_1),\cdots,2\arcsin(\sqrt\lambda\theta_d)\right), \quad \theta \in S^{d-1},
\label{S2xsqrtlambadomega}
\end{equation}
gives a parametric representation of
\begin{equation}
M_{\lambda} = \big\{x \in {\bf T}^d \, ; \, h(x) = \lambda \big\}.
\label{S4Mlambda}
\end{equation}
We equip $M_{\lambda}$ with the measure
\begin{gather*}
d\widetilde M_{\lambda} = \frac{(\sqrt\lambda)^{d-2}}{2}J(\sqrt{\lambda}\theta )d\theta , \\
J(y) = \chi(y){\mathop\prod_{j=1}^d}\,\frac{2}{\cos(x_j/2)} = \chi(y)
{\mathop\prod_{j=1}^d}\,\frac{2}{\sqrt{1-y_j^2}},
\end{gather*}
$\chi(y)$ being the characteristic function of $[-1,1]^d$.
Then we have
$$
dx = J(y)dy =  \,d\widetilde M_{\lambda} d\lambda , \quad
d\widetilde M_{\lambda}  = \frac{dM_{\lambda}}{|\nabla_x h(x)|},
$$
where $dM_{\lambda}$ is the measure on $M_{\lambda}$ induced from $dx$. Let $L^2(M_{\lambda})$ be the Hilbert space with inner product
$$
(\varphi,\psi)_{L^2(M_{\lambda})} = \int_{M_{\lambda}}
\varphi \,\overline{\psi}\,d\widetilde M_{\lambda}.
$$
We define $\mathcal F_0(\lambda)f = {\rm Tr}_{M_{\lambda}}f$,
 where ${\rm Tr}_{M_{\lambda}}$ is the trace on $M_{\lambda}$. More precisely,
\begin{equation}
\left(\mathcal F_0(\lambda)f\right)(\theta ) = f(x(\sqrt{\lambda}\theta )).
\label{S2F0define}
\end{equation}
It then follows for $R_0(z) = (H_0-z)^{-1}$
\begin{equation*}
\frac{1}{2\pi i}( ( R_0 (\lambda +i0 )- R_0 (\lambda -i0 )) f ,g )_{ L^2 ({\bf T} ^d) } = ( \mathcal{F}_0 (\lambda ) f , \mathcal{F}_0 (\lambda ) g)_{ L^2 (M_{\lambda } )} , 
\end{equation*}
for $ \lambda \in (0,d ) \setminus {\bf Z}$ and $ f,g \in C^1 ({\bf T}^d )$. We then have by (\ref{S2LAP})
\begin{equation}
\mathcal{F}_0 (\lambda ) \in {\bf B}(\mathcal B; L^2(M_{\lambda})).
\label{S2TraceOp}
\end{equation}
 Using this formula, we can derive the spectral representations of $H_0$ and $H$. However, we omit it.

%%%%%%%%%%%%%%%% subsection 3.3 %%%%%%%%%%%%%%%%%%

\subsection{Spectral representation on the lattice}
We define the distribution $ \delta ( h(x)-\lambda ) \in \mathcal{D}' ({\bf T}^d ) $ by 
$$
\int_{{\bf T}^d} f(x) \delta (h(x)-\lambda ) dx := \int_{ M_{\lambda } } f(x) \, d\widetilde{M} _{\lambda } , \quad f\in C^{\infty } ({\bf T}^d ) .
$$
Then, from the definition of $\mathcal{F}_0 (\lambda )^{\ast} $:
$$
\big( \mathcal{F}_0 (\lambda ) f , \, \phi \big) _{L^2 (M_{\lambda } )} = \big( f , \, \mathcal{F}_0 (\lambda )^{\ast} \phi \big) _{L^2 ({\bf T}^d )} , 
$$
we see that $\mathcal{F}_0 (\lambda )^{\ast} $ defines a distribution on ${\bf T}^d $ by the following formula 
$$
\mathcal{F}_0 (\lambda )^{\ast} \phi = \phi (x)\delta (h(x) -\lambda ) .
$$
Here the right-hand side makes sense when, for example, $\phi \in C^{\infty } (M_{\lambda } )$ and is extended to a $C^{\infty}$-function near $M_{\lambda }$.
Then $\widehat{\mathcal{F}}_0 (\lambda )^{\ast} \phi = \mathcal{U}^{\ast} \mathcal{F}_0 (\lambda )^{\ast} \phi $ is computed as 
\begin{gather}
\begin{split}
&(2\pi )^{-d/2} \int_{{\bf T}^d} e^{in \cdot x} \phi (x) \delta (h(x)-\lambda ) dx \\
=& (2\pi )^{-d/2} \int_{M_{\lambda } } e^{in \cdot x} \phi (x) \, d \widetilde{M}_{\lambda } \\
=& (2\pi )^{-d/2} \int_{S^{d-1} } e^{in \cdot x(\sqrt{\lambda } \theta )} \phi (x(\sqrt{\lambda } \theta )) \, \frac{( \sqrt{\lambda} )^{d-2} }{2} J(\sqrt{\lambda } \theta ) \, d\theta .
\end{split}
\label{S3psi0nlambdatheta2}
\end{gather}
In the lattice space, we define
$\widehat\psi^{(0)}(\lambda,\theta) = \big\{ \widehat\psi^{(0)}(n,\lambda,\theta) \big\}_{n\in{\bf Z}^d}$,
where 
\begin{equation}
\begin{split}
\widehat\psi^{(0)}(n,\lambda,\theta)  
&= (2\pi)^{-d/2}\frac{(\sqrt\lambda)^{d-2}}{2}e^{in\cdot x(\sqrt\lambda\theta)}J(\sqrt\lambda\theta) \\
&= (2\pi)^{-d/2}2^{d-1}(\sqrt\lambda)^{d-2}\chi(\sqrt\lambda\theta)
\frac{e^{in\cdot x(\sqrt{\lambda}\theta)}}{{\mathop\prod_{j=1}^d}\cos\big(x_j(\sqrt{\lambda}\theta)/2\big)}.
\end{split}
\label{S3psi0nlambdatheta}
\end{equation}
Here $\chi(y)$ is the characteristic function of $[-1,1]^d$, and $x(\sqrt{\lambda}\theta)$ is defined by (\ref{S2xsqrtlambadomega}).
By (\ref{S3psi0nlambdatheta2}) and (\ref{S3psi0nlambdatheta}), we have for $\phi \in L^2(M_{\lambda})$
\begin{equation}
\begin{split}
( \widehat{\mathcal F}_0(\lambda)^{\ast}\phi )(n) &=
(2\pi)^{-d/2}\int_{M_{\lambda}}e^{in\cdot x}\phi (x) \, d\widetilde M_{\lambda} \\
& =
\int_{S^{d-1}}{\widehat\psi}^{(0)}(n,\lambda,\theta ) \phi ( x (\sqrt{\lambda } \theta ) ) \, d\theta.
\end{split}
\nonumber
\end{equation}
We can also see for rapidly decreasing $\widehat f$ on ${\bf Z}^d$
\begin{equation}
( \widehat{\mathcal F}_0(\lambda)\widehat f )( x(\sqrt{\lambda } \theta ))=
(2\pi)^{-d/2}\sum_{n\in {\bf Z}^d} e^{-in\cdot x(\sqrt{\lambda}\theta )}\widehat f(n).
\nonumber
\end{equation}

The spectral representation for $\widehat H$ is constructed as follows. We put
\begin{equation}
\widehat{\mathcal F}^{(\pm)}(\lambda) = \widehat{\mathcal F}_0(\lambda)\left(1 - \widehat{V}\widehat{R}(\lambda \pm i0)\right), \quad
\lambda \in (0,d)\setminus {\bf Z}.
\label{S3Fpmlambdadefine}
\end{equation}
Then by (\ref{S2TraceOp}) and (\ref{S2LAP})
\begin{equation}
\widehat{\mathcal F}^{(\pm)}(\lambda) \in {\bf B}(\widehat{\mathcal B}\, ;L^2(M_{\lambda})).
\nonumber
\end{equation}
We define the operator $\widehat{\mathcal F}^{(\pm)}$  by $\big(\widehat{\mathcal F}^{(\pm)}f\big)(\lambda) = \widehat{\mathcal F}^{(\pm)}(\lambda)f$ for $f \in \widehat{\mathcal B}$.

%%%%%%%%%%%%%%%% Theorem 3.1 %%%%%%%%%%%%%%%%

\begin{theorem}
(1) $\widehat{\mathcal F}^{(\pm)}$ is uniquely extended to a partial isometry with initial set $\mathcal H_{ac}(\widehat H)$ and final set $L^2({\bf T}^d)$. Moreover it diagonalizes $\widehat H$:
\begin{equation}
\big(\widehat{\mathcal F}^{(\pm)}\widehat H\widehat f\big)(\lambda) = \lambda \big(\widehat{\mathcal F}^{(\pm)}\widehat f\big)(\lambda), \quad \widehat f \in \mathcal{H}_{ac} (\widehat H).
\end{equation}
(2) The  following inversion formula holds:
\begin{equation}
\widehat f = \mathop{\rm s-lim}_{N\to\infty}\int_{I_N}
\widehat{\mathcal F}^{(\pm)}(\lambda)^{\ast}\big(\widehat{\mathcal F}^{(\pm)}\widehat f\big)(\lambda)d\lambda, \quad \widehat f \in \mathcal H_{ac}(\widehat H),
\end{equation}
where $I_N$ is a union of compact intervals in $ (0,d)\setminus {\bf Z}$ such that $I_N \to (0,d)\setminus {\bf Z}$. \\
\noindent
(3) $\widehat{\mathcal F}^{(\pm)}(\lambda)^{\ast} \in {\bf B}(L^2(M_{\lambda})\,;\,\widehat{\mathcal B}^{\ast})$ is an eigenoperator for $\widehat H$ in the sense that
$$
(\widehat H - \lambda)\widehat{\mathcal F}^{(\pm)}(\lambda)^{\ast}\phi = 0, \quad \phi \in L^2(M_{\lambda}).
$$
(4) The wave operators
\begin{equation}
\widehat W^{(\pm)} = {\mathop{\rm s-lim}_{t\to\pm\infty}}\, e^{it\widehat H}e^{-it\widehat H_0}
\nonumber
\end{equation}
exist and are complete. Moreover,
\begin{equation}
{\widehat W}^{(\pm)} = \big(\widehat{\mathcal F}^{(\pm)}\big)^{\ast}\widehat{\mathcal F}_0.
\nonumber
\end{equation}
\end{theorem}

%%%%%%%%%%%%%%%%% subsection 3.4 %%%%%%%%%%%%%%%

\subsection{Scattering matrix}
 The scattering operator $\widehat S$ is defined by
\begin{equation}
\widehat S = \big(\widehat W_{+}\big)^{\ast}\widehat W_{-}.
\nonumber
\end{equation}
We conjugate it by the spectral representation. Let
\begin{equation}
\mathcal S = \widehat {\mathcal F}_0\widehat S(\widehat{\mathcal F}
_0)^{\ast},
\nonumber
\end{equation}
which is unitary on $L^2((0,d);L^2(M_{\lambda});d\lambda)$.  Since $\mathcal S$ commutes with $\widehat H_{0}$,
${\mathcal S}$ is written as a direct integral
\begin{equation}
{\mathcal S} = \int_{(0,d)}{\oplus}{\mathcal S}(\lambda)d\lambda.
\nonumber
\end{equation}
The S-matrix, ${\mathcal S}(\lambda)$, is unitary on $L^2(M_{\lambda})$ and has the
following representation.

%%%%%%%%%%%%%%%%%%% Theorem 3.2 %%%%%%%%%%%%%%

\begin{theorem} \label{Wholespacescatteringamp}
Let $\lambda \in (0,d)\setminus {\bf Z}$.  Then ${\mathcal S}(\lambda)$ is written as
\begin{equation}
{\mathcal S}(\lambda) = 1 - 2\pi i A(\lambda),
\nonumber
\end{equation}
where
\begin{equation}
A(\lambda) = \widehat{\mathcal F}_0(\lambda)\left(1 - \widehat V\widehat R(\lambda + i0)\right)\widehat V\widehat{\mathcal F}_0(\lambda)^{\ast} 
= \widehat{\mathcal F}^{(+)}(\lambda)\widehat V\widehat{\mathcal F}_0(\lambda)^{\ast} ,
\label{S3ScattAmpWholeSp}
\end{equation}
and is called the scattering amplitude.
\end{theorem}

%%%%%%%%%%%%%%%%%%%%% Section 4 %%%%%%%%%%%%%%%%%%%%%%

\section{Asymptotic expansion of the resolvent at infinity}

%%%%%%%%%%%%%% subsection 4.1 %%%%%%%%%%%%%%%%%%%%%%%

\subsection{Stationary phase method on a surface}
Let $S$ be a compact $C^{\infty}$-surface in ${\bf R}^d$ of 
codimension 1, and $dS$ the measure on $S$ induced from the Euclidean metric. For $a(x) \in C^{\infty}(S)$ and $k \in {\bf R}^d$, we put
\begin{equation}
I(k) = \int_Se^{ix\cdot k}a(x)dS.
\label{S4I(x)}
\end{equation}

%%%%%%%%%%%%%%%%%%%% Theorem 4.1 %%%%%%%%%%%%%%%%%%

\begin{theorem}
Let $N(x)$ be an outward unit normal field on $S$, and $W(x)$, 
$K(x)$ the Weingarten map and the Gaussian curvature at $x \in S$, respectively.
Assume that there exists a finite number of points 
$x^{(j)}_{\pm} \in S$, $j= 1, \cdots, {\nu}$, such that
$$
k/|k| = \pm N(x_{\pm}^{(j)}),
$$
and that $K(x_{\pm}^{(j)}) \neq 0$, $j = 1,\cdots,\nu$.
Then we have as $\rho = |k| \to \infty$
\begin{equation}
\begin{split}
I(k)  &=  \rho^{-(d-1)/2}\sum_{j=1}^{\nu}
e^{ik\cdot x_+^{(j)}}
A_{+}(x_+^{(j)})\\
& + \rho^{-(d-1)/2}\sum_{j=1}^{\nu}e^{ik\cdot x_-^{(j)}}
A_-(x_-^{(j)}) + O(\rho^{-(d+1)/2}),
\end{split}
\label{S4I(xi)expand}
\end{equation}
where
\begin{equation}
A_{\pm}(x) = (2\pi)^{(d-1)/2}|K(x)|^{-1/2}
e^{\mp {\rm sgn}\,W(x)\, \pi i/4}a(x).
\label{S4Apm(x)}
\end{equation}
and ${\rm sgn}\,W(x) = n_+ - n_-$, $n_+$ $(n_-)$ being the number of positive (negative) eigenvalues of $W(x)$.
\end{theorem}

For the proof, see Lemma 2.2 and appendix of \cite{Matsu68}.
See also \cite{Litt63}.
If $S$ is represented by $x_d = f(x ')$, $x '= (x_1,\cdots,x_{d-1})$, the Gaussian curvature is given by
\begin{equation}
K(x) = \Big(\sum_{i=1}^{d-1}\big( \frac{\partial f}{\partial x_i} (x') \big)^2 +1 \Big)^{-(d+1)/2}\det\Big(- \frac{\partial^2f}{\partial x_i\partial x_j}(x')\Big).
\label{S4Gausscurvature}
\end{equation}
For $d = 2$, the Gaussian curvature of the curve $f(x_1,x_2) = 0$ is computed as
\begin{equation}
\big|K(x_1,x_2)\big| = \frac{\big|f_{x_2x_2}\cdot f_{x_1}^2 - 2f_{x_1x_2}\cdot f_{x_1}f_{x_2} + f_{x_1x_1}\cdot f_{x_2}^2\big|}{(f_{x_1}^2 + f_{x_2}^2)^{3/2}}.
\label{S4Gausscurvatured=2}
\end{equation}

%%%%%%%%%%%%%%% subsection 4.2 %%%%%%%%%%%%%%%%%%%%

\subsection{Convexity of $M_{\lambda}$}
As will be seen below, the shape of $M_{\lambda}$ depends highly on the space dimension and $\lambda$.
We know that $\nabla h(x) \neq 0$ on $M_{\lambda}$ if $\lambda \not\in {\bf Z}$. Assume that at a point in $M_{\lambda}$,
$\partial h/\partial x_d = (\sin x_d)/2 \neq 0$. 
We take $x_1, \cdots, x_{d-1}$ as local coordinates, and differentiate $h(x) = \lambda$ to get
$$
\sin x_i + \sin x_d\frac{\partial x_d}{\partial x_i} = 0,
$$
$$
\delta_{ij}\cos x_j + \cos x_d \frac{\partial x_d}{\partial x_i}\frac{\partial x_d}{\partial x_j} + \sin x_d\frac{\partial^2 x_d}{\partial x_i\partial x_j} = 0,
$$
for $i, j = 1, \cdots, d-1$.
We put
$\varphi = \sum_{j=1}^dk_jx_j.$
Then we have on $M_{\lambda}$
$$
\frac{\partial\varphi}{\partial x_i} = k_i + k_d\frac{\partial x_d}{\partial x_i} = k_i - k_d\frac{\sin x_i}{\sin x_d},
$$
$$
\frac{\partial^2\varphi}{\partial x_i\partial x_j} = k_d\frac{\partial^2x_d}{\partial x_i\partial x_j} = - \frac{k_d }{(\sin x_d)^3}\left(\delta_{ij}\cos x_j(\sin x_d)^2 + \sin x_i\sin x_j\cos x_d\right).
$$
Suppose $\partial\varphi/\partial x_i = 0$, $i = 1, \cdots, d-1$. Then 
$$
k_i =\rho\sin x_i, \quad i = 1, \cdots, d,
$$
$$
\rho = |k|\left((\sin x_1)^2 + \cdots + (\sin x_d)^2\right)^{-1/2}.
$$
Therefore we have
\begin{equation}
\frac{\partial^2\varphi}{\partial x_i\partial x_j} = -
\frac{1}{(\sin x_d)^2 \rho }\left(\delta_{ij}k_d^2\cos x_j + 
k_ik_j\cos x_d\right).
\label{S4del2phidxidxj}
\end{equation}
 
Now let us compute the determinant $\det\left(\partial^2\varphi/\partial x_i\partial x_j\right)$.

\bigskip
\noindent
(1) {\it The case $d=2$}. Using $k_i = \rho\sin x_i$, we have
\begin{equation}
\begin{split}
k_2^2\cos x_1 + k_1^2\cos x_2 & =\rho^2 (\cos x_1 + \cos x_2)(1 - \cos x_1\cos x_2) \\
& = 2 \rho^2 (1-\lambda)(1 - \cos x_1\cos x_2).
\end{split}
\nonumber
\end{equation}
Since $\lambda \neq 1$, this vanishes if and only if $\cos x_1 = \cos x_2 = \pm 1$, i.e. $x_1 = 0$ or $\pi$, and $x_2 = 0$ or $\pi$. However in this case,
$h(x) = \sum_{i=1}^2\sin^2(x_i/2) \in {\bf Z}.$
This implies that 
\begin{equation}
\partial^2\varphi/\partial x_1^2 \neq 0 \quad {\rm for} \quad 
\lambda \in (0,1)\cup(1,2).
\end{equation}
Therefore $M_{\lambda}$ is a closed curve in ${\bf T}^2$, and convex in the fundamental domain ${\bf R}^2/(2\pi{\bf Z})^2$, as is seen from the figures (Figures 1, 2, 3) below.
Let us remark here, in view of Figure 3, in the case $ 1 < \lambda < 2 $, it is convenient to shift the fundamental domain so that $ {\bf R}^2 / (2 \pi {\bf Z} )^2= [0, 2 \pi ] ^2 $. 
\begin{figure}[hbtp]
\begin{minipage}{0.4\hsize}
  \begin{center}
  \includegraphics[width=4.5cm, bb=0 0 212 213]{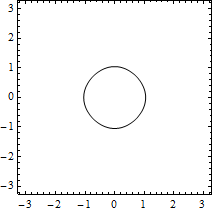}
  \end{center}
  \caption{$d=2$, $\lambda = 0.25$.}
\label{fig:math2}
 \end{minipage}
 \begin{minipage}{0.4\hsize}
  \begin{center}
   \includegraphics[width=4.5cm, bb=0 0 212 213]{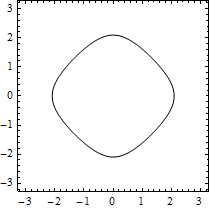}
  \end{center}
  \caption{$d=2$, $\lambda =0.75$.}
\label{fig:math3}
 \end{minipage}
\end{figure}

\begin{figure}[hbtp]
\centering
\includegraphics[width=4.5cm, bb=0 0 212 213]{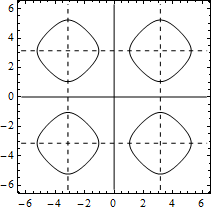}
\caption{$d=2$, $\lambda=1.25$.}
\label{fig:math5}
\end{figure}

\noindent
(2) {\it The case $d=3$}. By a direct computation, we have
\begin{equation}
\begin{split}
& \det\Big(\delta_{ij}k_3^2\cos x_j + 
k_ik_j\cos x_3\Big) \\
& = k_3^2\left(k_1^2\cos x_2\cos x_3 + k_2^2\cos x_3\cos x_1+ k_3^2\cos x_1\cos x_2\right),
\end{split}
\nonumber
\end{equation}
which can vanish when e.g. $\cos x_1 = \cos x_2 = 0$, $\cos x_3 = 1/2$.
Therefore in 3-dimensions, $M_{\lambda}$ may not be convex.
The following Figures 4, 5, 6 explain the situation in 3-dimensions.

\begin{figure}[hbtp]
\begin{minipage}{0.4\hsize}
\centering
\includegraphics[width=4.5cm, bb=0 0 187 209]{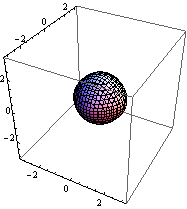}
\caption{$d=3, \lambda = 0.45$.}
\label{fig:math6}
\end{minipage}
\begin{minipage}{0.4\hsize}
\centering
\includegraphics[width=4.5cm, bb=0 0 180 220]{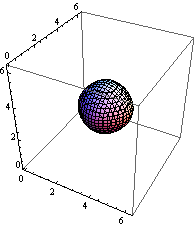}
\caption{$d=3, \lambda=2.55$.}
\label{fig:math7}
\end{minipage}
\end{figure}
\begin{figure}[hbtp]
\centering
\includegraphics[width=4.5cm, bb=0 0 187 209]{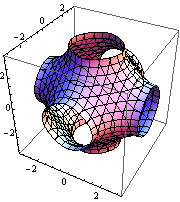}
\caption{$d=3, \lambda=1.45$.}
\label{fig:math8}
\end{figure}

Here, we note the following simple lemma.

%%%%%%%%%%%%%%%%%%% Lemma 4.2 %%%%%%%%%%%%%%%%%%

\begin{lemma}
If $-1 \leq y_i \leq 1$, $i = 1,\cdots,d$, and $d-1 < y_1+ \cdots + y_d < d$, we have $y_i > 0$, $ i = 1,\cdots d$.
\end{lemma}

Proof. Suppose e.g. $y_d \leq 0$. Then
$$
y_1 + y_2 + \cdots + y_d \leq y_1+ \cdots + y_{d-1} \leq d-1,
$$
which is a contradiction. \qed

\medskip
By (\ref{S4del2phidxidxj}), we have
$$
\sum_{i,j=1}^{d-1} \frac{\partial^2\varphi}{\partial x_i\partial x_j}\xi_i\xi_j = 
-\frac{1}{(\sin x_d )^2 \rho }\left(k_d^2\sum_{i=1}^{d-1}  \big( \cos x_i \big) \xi_i^2 + 
\big( \cos x_d \big)  \Big( \sum_{i=1}^{d-1} k_i\xi_i \Big) ^2 \right),
$$
which has a definite sign if $\cos x_i > 0$, $i = 1,\cdots,d$ and $\sin x_d >0$.
By virtue of Lemma 4.2, it happens for $0< \lambda < 1/2$.
Let us also note that for $ d-1/2 < \lambda < d$, we have the same conclusion since $ \cos x_i <0 \ (i=1, \cdots , d )$, $ \sin x_d <0$.
Recall that when $ d\geq 3$ the definition of the Gaussian curvature depends on the choice of direction of the unit normal $N(x) $ on $S$.
We choose $N(x)$ in such a way that $K(x)>0$ on $S$.

With this convention, we have proven the following lemma. Recall the interval $I_d$ defined by (\ref{S1IntervalI}).

%%%%%%%%%%%%%%%%%%%%%% Lemma 4.3 %%%%%%%%%%%%%%%%%%%%%%%

\begin{lemma} 
If $\lambda \in I_d$, all the principal curvature of $M_{\lambda}$ are positive.
\end{lemma}

As has been noted above, in the case $ 1<\lambda <2 \ (d=2 )$ or $ d-1/2 <\lambda <d \ (d\geq 3 )$, we should shift the fundamental domain so that $ {\bf R}^d /(2\pi {\bf Z})^d= [0 , 2\pi ]^d $ (See Figures 3, 4, 5).
To fix the idea, in the sequel, we deal with the case ${\bf R}^d / (2\pi {\bf Z})^d = {\bf T}^d = [-\pi , \pi ]^d$.

Under the assumption of Lemma 4.3, $M_{\lambda}$ is strictly convex. Let $N(x)$ be the unit normal field on $M_{\lambda}$ specified as above. Then for any $\omega \in S^{d-1}$, there exists a unique pair of points
 $x_{\pm}(\lambda,\omega)$ in $M_{\lambda}$ such that 
\begin{equation}
N(x_{\pm}(\lambda,\omega)) = \pm \omega.
\label{S4N(x)=k/|k|}
\end{equation}
Since $N(-x) = - N(x)$, we see that
$x_-(\lambda,\omega) = - x_+(\lambda,\omega)$. Therefore, we let
\begin{equation}
x_{\pm}(\lambda,\omega) = \pm x_{\infty }(\lambda,\omega).
\label{S4xpm=pmx}
\end{equation}

We can now compute the asymptotic expansion of the free resolvent
\begin{gather}
\big(\widehat R_0(z)\widehat f\big)(m) = \sum_{n\in{\bf Z}^d}r_0(m-n,z)\widehat f(n), \\
r_0(k,z) = (2\pi)^{-d}\int_{{\bf T}^d}\frac{e^{ik\cdot x}}{h(x) - z}dx.
\end{gather}
We put
\begin{equation}
\omega_k = k/|k|, \quad k \in {\bf R}^{d}\setminus\{0\}.
\label{S4omegak}
\end{equation}

%%%%%%%%%%%%%%%%%% Lemma 4.4 %%%%%%%%%%%%%%%%%%%%%

\begin{lemma}
Assume $\lambda \in I_d$. 
Then we have as $|k| \to \infty$
\begin{equation}
\begin{split}
& r_0(k,\lambda\pm i0) \\
= & \pm i(2\pi|k|)^{-(d-1)/2}
e^{\pm i(k\cdot x_{\infty} (\lambda,\omega_k) - (d-1)\pi/4)}\frac{K(x_{\pm}(\lambda,\omega_k))^{-1/2}}{|\nabla_x h(x_{\pm}(\lambda,\omega_k))|} + O(|k|^{-(d+1)/2}).
\end{split}
\nonumber
\end{equation}
\end{lemma}

Proof. Take $\epsilon > 0$ small enough so that
\begin{equation}
(\lambda - 2\epsilon,\lambda + 2\epsilon) \subset 
\left\{
\begin{split}
& (0,1) , \quad d = 2, \\
& (0,1/2), \quad d \geq 3.
\end{split}
\right.
\nonumber
\end{equation}
Let $\chi(t) \in C_0^{\infty}({\bf R})$ be such that $\chi(t) = 1$ for $|t|< \epsilon /2 $, $\chi(t) = 0$ for $|t| > \epsilon$, and assume that $|{\rm Re}\,z - \lambda| < \epsilon /4$. We split $r_0(k,z)$ into two parts
\begin{equation}
r_0(k,z) = A(k,z) + B(k,z),
\nonumber
\end{equation}
\begin{equation}
A(k,z) = (2\pi)^{-d}\int_{{\bf T}^d}
\frac{\chi(h(x)-\lambda)}{h(x)-z}e^{ik\cdot x}dx.
\nonumber
\end{equation}
Then, by integration by parts, for all $N > 0$
\begin{equation}
B(k,z) = O(|k|^{-N}), \quad |k| \to \infty.
\nonumber
\end{equation}
Letting $S(t) = \big\{x \in {\bf T}^d\, ; \, h(x) = t\big\}$, we write $A(k,z)$  as
\begin{equation*}
A(k,z) = (2\pi)^{-d}\int_{\lambda -  \epsilon }^{\lambda +  \epsilon }\frac{a(t,k)}{t-z}dt, \quad
a(t,k) = \int_{S(t)}e^{ik\cdot x}\frac{\chi(t-\lambda)}{|\nabla_x h(x)|}dS(t).
\end{equation*}
We then have
\begin{equation}
\int_{ \lambda -  \epsilon }^{ \lambda + \epsilon }\frac{a(t,k)}{t - \lambda \mp i0} dt = \pm i \pi a(\lambda,k) + {\rm p.v.}\int_{\lambda - \epsilon }^{ \lambda + \epsilon}
\frac{a(t,k)}{t - \lambda}dt.
\label{S4limpmi0int}
\end{equation}
By Theorem 4.1, for $t \in (\lambda-\epsilon,\lambda+ \epsilon)$, $a(t,k)$ admits the asymptotic expansion
\begin{equation}
a(t,k) = a_0(t,k) + O(|k|^{-(d+1)/2}),
\nonumber
\end{equation}
\begin{equation}
\begin{split}
a_0(t,k) &= \Big(\frac{2\pi}{|k|}\Big)^{(d-1)/2}
e^{ik\cdot x_{\infty} (t,\omega_k) - (d-1)\pi i/4}
\chi (t-\lambda ) \frac{K(x_+(t,\omega_k))^{-1/2}}{|\nabla_x h(x_+(t,\omega_k))|} \\
& + \Big(\frac{2\pi}{|k|}\Big)^{(d-1)/2}
e^{-ik\cdot x_{\infty } (t,\omega_k) + (d-1)\pi i/4}
\chi (t- \lambda ) \frac{ K(x_-(t,\omega_k))^{-1/2}}{|\nabla_x h(x_-(t,\omega_k))|}\\
&=: a_0^{(+)}(t,k) + a_0^{(-)}(t,k),
\end{split}
\label{S4a0pm(t,k)}
\end{equation}
where $x_{\pm}(t,\omega_k)$ is a stationary phase point on $S(t)$. 

We compute the asymptotic expansion of the 2nd term of the right-hand side of (\ref{S4limpmi0int}). Differentiating $h(x_{\pm}(t,\omega_k)) = t$, we have
$$
\nabla_x h(x_{\pm}(t,\omega_k))\cdot\partial_t x_{\pm}(t,\omega_k) = 1.
$$
Therefore, letting 
$$
s = \omega_k\cdot x_{\pm}(t,\omega_k) - \omega_k\cdot x_{\pm}(\lambda,\omega_k),
$$
we have
$$
\frac{ds}{dt} = \omega_k\cdot \partial_t x_{\pm}(t,\omega_k) = \frac{\nabla_x h(x_{\pm}(t,\omega_k))}{|\nabla_x h(x_{\pm}(t,\omega_k)|}\cdot\partial_tx_{\pm}(t,\omega_k) = \frac{1}{|\nabla_x h(x_{\pm}(t,\omega_k))|},
$$
which implies
$$
t - \lambda = s|\nabla_x h(x_{\pm}(\lambda,\omega_k))| + O(s^2).
$$
We then have 
$$
\frac{1}{t-\lambda}\frac{\chi (t-\lambda ) K(x_{\pm}(t,\omega_k))^{-1/2}}{|\nabla_x h(x_{\pm}(t,\omega_k)|}\frac{dt}{ds} = 
\frac{b_{\pm}(s,\omega_k)}{s},
$$
where $b_{\pm}(s,\omega_k)$ is a smooth function such that
$$
b_{\pm}(0,\omega_k) = \frac{K(x_{\pm}(\lambda,\omega_k))^{-1/2}}{|\nabla_x h(x_{\pm}(\lambda,\omega_k))|}.
$$
Taking $\delta > 0$ small enough, we have by integration by parts
\begin{equation}
\begin{split}
{\rm p.v.}\int_{-\delta}^{\delta}\frac{e^{\pm i|k|s}}{s}b_{\pm}(s,\omega_k)ds 
&= \pm 2i \int_0^{|k|\delta}\frac{\sin s}{s}ds\, b_{\pm}(0,\omega_k) + O(|k|^{-1})\\
&= \pm \pi i\, b_{\pm}(0,\omega_k) + O(|k|^{-1}),
\end{split}
\nonumber
\end{equation}
which implies
\begin{equation}
\begin{split}
& {\rm p.v.}\int_{\lambda - \epsilon}^{\lambda + \epsilon}\frac{a_0^{(\pm)}(t,k)}{t-\lambda}dt \\
& =  \Big(\frac{2\pi}{|k|}\Big)^{(d-1)/2}
e^{\pm i k\cdot x_{\infty }(\lambda,\omega_k) \mp (d-1)i\pi /4}{\rm p.v.}
\int_{-\delta}^{\delta}\frac{e^{\pm i|k|s}}{s}b_{\pm}(s,\omega_k)ds + O(|k|^{-(d+1)/2}) \\
& =  \pm i\pi\Big(\frac{2\pi}{|k|}\Big)^{(d-1)/2}
 e^{\pm i k\cdot x_{\infty } (\lambda,\omega_k) \mp (d-1)i\pi /4}\frac{K(x_{\pm}(\lambda,\omega_k))^{1/2}}{|\nabla_x h(x_{\pm}(\lambda,\omega_k))|} + O(|k|^{-(d+1)/2}) .
\end{split}
\label{S4iintpvapmtk}
\end{equation}
Plugging (\ref{S4limpmi0int}), (\ref{S4a0pm(t,k)}) and (\ref{S4iintpvapmtk}), we obtain the lemma. \qed

%%%%%%%%%%%%%%%%% Lemma 4.5 %%%%%%%%%%%%%%%%%

\begin{lemma}
We have as $|m| \to \infty$
$$
(m-n)\cdot x_{\pm}(\lambda,\omega_{m-n}) = (m-n)\cdot x_{\pm}(\lambda,\omega_{m}) + O(|m|^{-1}).
$$
\end{lemma}

Proof.  We extend $x_{\pm}(\lambda,k)$ as a function of homogeneous degree 0 in $k$. Letting $\epsilon = 1/|m|$, we have 
$$
\omega_{m-n} = (\omega_m - \epsilon n)/|\omega_m - \epsilon n| = \omega_m + \epsilon((\omega_m\cdot n)\omega_m - n) + O(\epsilon^2).
$$
Using $h (x_{\pm}(\lambda,\omega_{m-n})) = \lambda$, we have
$$
\nabla_x h(x_{\pm}(\lambda,\omega_{m-n}))\cdot\frac{d}{d\epsilon}x_{\pm}(\lambda,\omega_{m-n})\Big|_{\epsilon=0} = 0.
$$
Since $\nabla_x h (x_{\pm}(\lambda,\omega))$ is parallel to $\omega$, we then have
$$
\omega_m\cdot\frac{d}{d\epsilon}x_{\pm}(\lambda,\omega_{m-n})\Big|_{\epsilon=0} = 0,
$$
which implies
$$
m\cdot x_{\pm}(\lambda,\omega_{m-n}) = m\cdot x_{\pm}(\lambda,\omega_m) + O(
|m|^{-1}),
$$
and the lemma follows immediately. \qed

\bigskip
Lemmas 4.4 and 4.5 imply the following lemma.

%%%%%%%%%%%%%%% Lemma 4.6 %%%%%%%%%%

\begin{lemma} \label{R0lambdakinfty}
If $\lambda \in I_d$ and $\widehat f(n)$ is compactly supported, we have as $|k| \to \infty$
\begin{equation}
\begin{split}
 & \left(\widehat R_0(\lambda \pm i0)\widehat f\right)(k) \\
 & =  e^{\pm (3-d)\pi i/4}
(2\pi|k|)^{-(d-1)/2}
e^{\pm ik\cdot x_{\infty }(\lambda,\omega_k)}
a_{\pm}(\lambda,\omega_k)
\sum_ne^{\mp in\cdot x_{\infty } (\lambda,\omega_k)}\widehat f(n) \\
&  + O(|k|^{-(d+1)/2}),
\end{split}
\nonumber
\end{equation}
\begin{equation}
a_{\pm}(\lambda,\omega_k) = \frac{K(x_{\pm}(\lambda,\omega_k))^{-1/2}}{|\nabla_x h(x_{\pm}(\lambda,\omega_k))|}.
\label{S4apm}
\end{equation}\label{asympt_r0}
\end{lemma}

 Recalling the definition of $x(\sqrt{\lambda}\theta )$ in 
(\ref{S2xsqrtlambadomega}) and the fact that the Gauss map is a diffeomorphism for a strictly convex surface, define $\theta(\lambda,\omega)$ by the relation 
$x(\sqrt{\lambda}\,\theta(\lambda,\omega)) = x_{\infty }(\lambda,\omega)$, i.e.
 \begin{equation}
\theta_j (\lambda,\omega) = \frac{1}{\sqrt{\lambda}}\sin\Big(\frac{1}{2}\,x_{\infty j}(\lambda,\omega)\Big) ,\quad j=1, \cdots ,d.
\label{S4defineomegapmj}
\end{equation}
We define the reparametrized Fourier transforms $\widehat{\mathcal G}_0(\lambda)$ and  $\widehat{\mathcal G}^{(\pm)}(\lambda)$ by
\begin{equation}
\left(\widehat{\mathcal G}_{0}(\lambda) \widehat{f} \right)(\omega) = \left(
\widehat{\mathcal F}_0(\lambda)\widehat f\right)(\theta(\lambda,\omega)),
\label{S4G0pmjlambdadefine}
\end{equation}
\begin{equation}
\widehat{\mathcal G}^{(\pm)}(\lambda) = \widehat{\mathcal G}_{0}(\lambda)(1 - \widehat V\widehat R(\lambda \pm i0)).
\label{S4Gpmjlambdadefine}
\end{equation}
Lemma 4.6, the definition (\ref{S3Fpmlambdadefine}) and the resolvent equation imply the following theorem.

%%%%%%%%%%%%%%% Theorem 4.7 %%%%%%%%%%

\begin{theorem} \label{Rlambdapmi0Asympto}
If $\lambda \in I_d$ and $\widehat f(n)$ is compactly supported, we have as $|k| \to \infty$
\begin{equation}
\begin{split}
& \left(\widehat R(\lambda \pm i0)\widehat f\right)(k) \\
 = & \ e^{\pm (3-d)\pi i/4}\sqrt{2\pi}|k|^{-(d-1)/2}
e^{\pm ik\cdot x_{\infty } (\lambda,\omega_k)}
a_{\pm}(\lambda,\omega_k)
\left(\widehat{\mathcal G}^{(\pm)} (\lambda )\widehat{f} \right) (\pm \omega_k ) \\
&  + O(|k|^{-(d+1)/2}).
\end{split}
\nonumber
\end{equation}
\label{S4expansion_resol2}
\end{theorem}

%%%%%%%%%%%%%%%%%%%%%%%%% Section 5 %%%%%%%%%%%%%%%%%%%%%%%%

\section{Radiation conditions on ${\bf Z}^d$}

The aim of this section is to introduce the radiation condition (Definition \ref{RadCond}) and prove the uniqueness theorem (Theorem \ref{RadCondUnique}).

%%%%%%%%%%%%%%%%%%%%%%%% subsection 5.1 %%%%%%%%%%%%%%%%%

\subsection{Green's formula}
For $m,n \in {\bf Z}^d $, we write $m\sim n $, if $|m-n|=1 $, i.e. there exists $j$ such that $m= n \pm e_j$.
We define the discrete Laplacian $\Delta_{disc} $ on ${\bf Z}^d $ by 
\begin{equation}
(\Delta_{disc} \widehat{u})(n) = - (\widehat H_0\widehat u)(n) =\frac{1}{4} \sum_{ m \sim n} \big(\widehat{u} (m)-\widehat{u} (n) \big) . \label{discrete_laplacian} 
\end{equation}
A set $D \subset {\bf Z}^d $ is said to be {\it connected} if for any $m,n\in D$, there exist $m^{(j)} \in D $, $j=0, \cdots , k $ such that $m^{(j)} \sim m^{(j+1 )} $, $j=0 , \cdots , k-1 $, and $m^{(0)} =m$, $m^{(k) }=n$.
A connected subset $D\subset {\bf Z}^d$ is called a {\it domain}.
For a domain $D\subset {\bf Z}^d$, we define
\begin{gather}
{\rm deg}\,(n) = \, ^{\#} \big\{m \in D\, ; \,m \sim n \big\}, \quad n \in D,
\label{S5deg(n)}
\\
\stackrel{\circ}D \, = \big\{n \in D\, ; \, {\rm deg}\,(n) = 2 d \big\},
\label{S5Interior}
\\
\partial D = \big\{ n \in D\, ;\, {\rm deg}\,(n) < 2 d \big\}.
\label{S5partialD}
\end{gather}
The normal derivative at the boundary is defined by
\begin{equation}
\left( \partial_{\nu}^D \widehat{u} \right) (n) =\frac{1}{4} \sum_{m\in \stackrel{\circ}D, m\sim n  } \big( \widehat{u} (n)-\widehat{u}(m) \big) , \quad n\in \partial D . \label{normal_deri} 
\end{equation}
Note that, compared with (\ref{discrete_laplacian}), $m $ and $n$ are interchanged.
Then the following Green's formula holds (see e.g \cite{Du84} and \cite{IsMo}): 
\begin{equation} 
\begin{split}
&\sum_{n\in \stackrel{\circ}D} \big( (\Delta_{disc} \widehat{u} )(n) \cdot \widehat{v} (n) -\widehat{u}(n) \cdot (\Delta_{disc} \widehat{u} )(n) \big) \\
&=\sum_{ n\in \partial D} \big( (\partial^D_{\nu } \widehat{u} )(n) \cdot \widehat{v} (n)-\widehat{u} (n) \cdot (\partial^D_{\nu} \widehat{v})(n) \big) . 
\label{green} 
\end{split} 
\end{equation}

%%%%%%%%%%%%%%%%%%%%%%%%% subsection 5.2 %%%%%%%%%%%%%%5

\subsection{Radiation condition}

%
%%%%%%%%%%%%%%%%%%%%%%%%%%
For $m,n$ such that $m\sim n$, we define the difference operator $\partial _{m-n} $ by \begin{equation*}
\big( \partial_{m-n} \widehat{f} \big) (n)= \widehat{f} (m)-\widehat{f}(n). \end{equation*}
%
%%%%%%%%%%%%%%%%%%%%% Lemma 5.1 %%%%%%%%%%%%%%%%%
%
\begin{lemma}
(1) Let $n(s)=n+s(m-n)$, where $m\sim n$.
Then we have 
\begin{equation*}
\partial_{m-n } (n\cdot x_{\infty} (\lambda , \omega_n ) )=\int_0^1 (m-n )\cdot x_{\infty } (\lambda , \omega_{n(s)}) ds. 
\end{equation*}
(2) If $m\sim n$, we have as $|n|\rightarrow \infty $ 
\begin{gather*}
\partial_{m-n} (n\cdot x_{\infty}(\lambda , \omega_n ) ) =(m-n )\cdot x_{\infty } (\lambda , \omega_n )+O(|n|^{-1} ) , \\
\partial_{m-n} \left( e^{in\cdot x_{\infty } (\lambda , \omega_n )} \right)=\left( e^{i(m-n) \cdot x_{\infty} (\lambda , \omega_n )} -1 \right) e^{in\cdot  x_{\infty } (\lambda , \omega_n )} +O(|n|^{-1} ) . 
\end{gather*}
\label{section5_estimate1} 
\end{lemma}

Proof. Differentiating $h(x_{\infty} (\lambda , \omega_{n(s) } ) )=\lambda $, we have 
\begin{equation*}
(\nabla _x h)(x_{\infty } (\lambda , \omega_{n (s)} )) \cdot \frac{d}{ds} x_{\infty } (\lambda , \omega_{n(s)} ) =0. 
\end{equation*}
Since $ (\nabla_x h)( x_{\infty } (\lambda , \omega_{n(s)} ) ) $ is parallel to $n(s)$, we then have 
\begin{equation*}
n(s) \cdot \frac{d}{ds} x_{\infty } (\lambda , \omega_{n(s)} ) =0, 
\end{equation*}
which implies 
\begin{equation*}
\frac{d}{ds} (n(s) \cdot x_{\infty } (\lambda , \omega_{n(s) } )) = (m-n )\cdot x_{\infty } (\lambda , \omega_{n(s)} ) . 
\end{equation*}
Integrating this equality, we obtain (1).
Since $ \omega_{n(s)} =\omega_n +O(|n|^{-1} )$, (2) follows from (1). \qed

\medskip

We now introduce the rectangular domain $D(R) $ such that
\begin{equation}
\stackrel{\circ}{D(R)}=\big\{ n\in {\bf Z}^d \ ; \ n \in [-R , R]^d \big\}, \quad R>0,
\label{RectAngDom}
\end{equation}
and the {\it radial derivative} $\partial_{rad} $ by 
\begin{gather}
(\partial _{rad}\, \widehat{u} )(k)=\frac{1}{4} \sum_{m\in \partial D(R(k) ), m\sim k  } (\widehat{u} (m)-\widehat{u}(k) ),
\label{radial_deri} \\
R(k)=  \max_{1\leq j\leq d}\, |k_j|, \quad k\in {\bf Z}^d.
\label{R(k)define}
\end{gather}
We put
\begin{equation}
A_{\pm} (\lambda , \omega_k )= \frac{1}{4} \sum_{m \in \partial D(R(k)), m\sim k   } \left( e^{ \pm i (m-k) \cdot x_{\infty } (\lambda , \omega_k ) } -1 \right) , \quad 
\omega_k = \frac{k}{|k|}.
\label{rad_coeff} 
\end{equation}

%%%%%%%%%%%%%%%%%%% Lemma 5.2 %%%%%%%%%%%%%%%%%

\begin{lemma}
(1) The right-hand side of (\ref{rad_coeff}) does not depend on $|k|$.\\
(2) There exists a constant $ \epsilon _0 (\lambda )>0 $ such that \begin{equation*}
\pm \mathrm{Im} A_{\pm} (\lambda , \omega_k ) > \epsilon _0 (\lambda ) , \end{equation*}
for any $\omega_k$.
\label{section5_estimate2} 

\end{lemma}
Proof. If $ m\sim k $, $m\in \partial D(R(k)) $, then $m-k =\pm e_j $ for some $j$. This $\pm e_j$ depends only on $\omega_k$, which proves (1).

Recall that $\nabla h (x)=\frac{1}{2} (\sin x_1 , \cdots , \sin x_d )$, hence letting $\omega_{k,j} $ be the $j$-th component of $\omega _k $, we have \begin{equation*}
\sin (x_{\infty j} (\lambda , \omega_k )) = c \omega _{k,j} \end{equation*}
for some constant $c>0$.
Suppose $m\sim k $, $m \in \partial D(R(k)) $.
If $\omega_{k,j} >0$, then either $m_j=k_j $ or $m_j= k_j +1 $.
If $\omega_{k,j} <0 $, then either $m_j=k_j $ or $m_j =k_j -1 $.
We then have that $\sin ((m-k)\cdot x_{\infty } (\lambda , \omega_k )) =c|\omega_{k,j}| $ for some $j$ such that $\omega_{k,j} \not=0 $.
Since 
\begin{gather*} 
\begin{split}
\pm \mathrm{Im} A_{\pm } (\lambda , \omega_{k} ) &=\frac{1}{4} \sum_{  m \in \partial D(R(k)), m\sim k } \sin ( (m-k)\cdot x_{\infty } (\lambda , \omega_k )) \\
&=\frac{c}{4} \sum_{m\in \partial D(R(k)), m\sim k } |\omega_{k,j} |, 
\end{split} 
\end{gather*}
and $\sum_j \omega_{k,j}^2 =1 $, the lemma follows. \qed

\medskip

Let us introduce two auxiliary norms, $\widehat{\mathcal{B}}^{\ast} _{{\bf R} } $-norm and $\widehat{\mathcal{B}}^{\ast} _{{\bf Z} }$-norm, on $\widehat{\mathcal{B}}^{\ast}$ by
\begin{gather*}
\| \widehat{u} \|^2 _{\widehat{\mathcal{B}}^{\ast}_{{\bf R}} }= \sup_{R>1,R\in {\bf R}} \frac{1}{R} \sum_{n\in \stackrel{\circ}{D(R)} } |\widehat{u} (n)|^2 , \\
\| \widehat{u} \|^2_{\widehat{\mathcal{B}}^{\ast} _{{\bf Z}}} =\sup_{\rho >1, \rho \in {\bf Z} }\frac{1}{\rho } \sum_{ n \in \stackrel{\circ}{D(\rho )} } |\widehat{u} (n) |^2 . 
\end{gather*}

%%%%%%%%%%% Lemma 5.3 %%%%%%%%%%%%%%%%%%%%

\begin{lemma}
 These three norms $\|  \cdot  \|_{\widehat{\mathcal{B}}^{\ast} } $, $\|  \cdot \|_{\widehat{\mathcal{B}}^{\ast} _{{\bf R}}} $,  and $\| \cdot \| _{\widehat{\mathcal{B}}^{\ast} _{{\bf Z} } }$ are equivalent.
\label{normeq2}
\end{lemma}

Proof. Let $A(R) = \{z\in{\bf C}^d;(\sum_{j=1}^d|z_j|^2)^{1/2} <R\}$, $B(R) = \{z\in{\bf C}^d;\max_{j}|z_j|<R\}$. Then there is a constant $\delta >0$ such that $A(\delta R) \subset B(R) \subset A(R/\delta)$, $\forall R > 0$. This implies
\begin{equation}
\frac{1}{R}\sum_{|n|<\delta R}|\widehat u(n)|^2 \leq 
\frac{1}{R}\sum_{n\in \stackrel{\circ}{D(R)}}|\widehat u(n)|^2 \leq
\frac{1}{R}\sum_{|n|<R/\delta}|\widehat u(n)|^2.
\nonumber
\end{equation}
Taking the supremum with respect to $R > \delta$ or $R> 1/\delta$, we get the equivalence of $\|\cdot  \|_{\widehat{\mathcal{B}}^{\ast} }$ norm and  $\|  \cdot \|_{\widehat{\mathcal{B}}^{\ast} _{{\bf R}}}$ norm.

Next we show the equivalence of  $\|\cdot \| _{\widehat{\mathcal{B}}^{\ast} _{{\bf R} } }$ norm and $\| \cdot \| _{\widehat{\mathcal{B}}^{\ast} _{{\bf Z} }}$ norm. Note that $f(r) = \sum_{n \in \stackrel{\circ}{D(r)}}|\widehat u(n)|^2$ is a right-continuous non-decreasing step function on $(0,\infty)$ with jump at integers.
For $R>1$, we take $\rho (R)= [R]$ = the largest positive integer such that $\rho (R) \leq R$. Then we  have
\begin{equation*}
\sup_{R>1} \frac{1}{R} \sum_{ n\in \stackrel{\circ}{D(R)}} |\widehat{u} (n)|^2 \leq \sup_{R>1} \frac{1}{\rho (R)} \sum_{n \in \stackrel{\circ}{D(\rho (R) )}} |\widehat{u} (n)|^2 .
 \end{equation*}
The converse inequality is proven by the following inequality
\begin{equation*}
\sup_{R>1} \frac{1}{\rho(R)} \sum_{ n \in \stackrel{\circ}{D(\rho (R) )}} |\widehat{u} (n)|^2 \leq  \sup_{R>1} \frac{2}{R} \sum_{n\in \stackrel{\circ}{D(R)}} |\widehat{u} (n)|^2. 
\qed
\end{equation*}

%%%%%%%%%%%%%%%% Lemma 5.4 %%%%%%%%%%%%%%%%%

\begin{lemma}
(1) \ If $\widehat{f} \in \ell^{\infty } ({\bf Z}^d )$ satisfies $|\widehat{f} (n)| \leq C(1 +|n| )^{-(d-1 )/2 } $, then 
\begin{equation}
\sup_{R>1} \frac{1}{R} \sum_{|n|<R} |\widehat{f} (n)|^2 <\infty , \quad \text{i.e.} \quad \widehat{f} \in \widehat{\mathcal{B}}^{\ast} . \label{besov_ineq1}
\end{equation}
(2) \ If $|\widehat{f} (n)|\leq C(1+|n|)^{-(d-1)/2 -\epsilon }$, $\epsilon >0$, then 
\begin{equation}
\lim_{R \to \infty } \frac{1}{R} \sum_{|n| <R} |\widehat{f} (n)|^2 =0 . \label{besov_ineq2} 
\end{equation}
\label{section5_besov_asymptotic}
\end{lemma}
Proof.
 We compute the norm $\|\widehat f \|_{\widehat{\mathcal{B}}^{\ast} _{{\bf Z}} } $. We first show 
\begin{equation}
\sum_{n\in \stackrel{\circ}{D(\rho)} \setminus \stackrel{\circ}{D(\rho -1 )}} |\widehat{f} (n)|^2 = O(1),
\label{lemma5_5_2} 
\end{equation}
as $\rho \to \infty$. In fact, for any $\rho \in {\bf Z}, \, \rho > 1$ and $n\in \, \stackrel{\circ}{D(\rho )} \setminus \stackrel{\circ}{D(\rho -1 )} $, we have
$\rho -1 < |n| \leq \sqrt{d} \rho $.  
Since $^{\#} \big\{n \in \, \stackrel{\circ}{D(\rho) } \setminus \stackrel{\circ}{D(\rho-1)} \big\} =  (2\rho +1 )^d -(2\rho -1)^d \leq C\rho^{d-1}$,
\begin{equation}
\sum_{ n \in \stackrel{\circ}{D(\rho )} \setminus \stackrel{\circ}{D(\rho -1 )}} |\widehat{f} (n)|^2 \leq C\rho^{ -(d-1 )}  \, ^{\#} \big\{n \in \, \stackrel{\circ}{D(\rho) } \setminus \stackrel{\circ}{D(\rho-1)} \big\} \leq C. 
\nonumber
\end{equation}
On the other hand, since
\begin{equation*}
\sum_{n\in \stackrel{\circ}{D(R )}} |\widehat{f} (n)|^2 = \sum_{\rho =1}^R \sum_{ n \in \stackrel{\circ}{D(\rho )} \setminus \stackrel{\circ}{D(\rho -1 )}} |\widehat{f} (n)|^2 +|\widehat{f} (0)|^2 , 
\end{equation*}
for every positive integer $R$, we have
$\sum_{n\in \stackrel{\circ}{D(R)}} |\widehat{f} (n)|^2 =O(R)$
by (\ref{lemma5_5_2}).
This proves (1) by Lemma \ref{normeq2}.

Assume $|\widehat{f} (n) |\leq C(1+|n|)^{-(d-1)/2 -\epsilon } $ for some $\epsilon >0$. By the similar computation, we have 
$\sum_{n\in \stackrel{\circ}{D(R)} } |\widehat{f} (n)|^2 =o(R)$, which proves (2). \qed

\medskip

For $\widehat{f},\ \widehat{g} \in \widehat{\mathcal{B}}^{\ast} $, we write 
\begin{equation}
\widehat{f} \simeq \widehat{g} \Longleftrightarrow \lim_{R \to \infty} \frac{1}{R} \sum_{|n|<R} |\widehat{f} (n)-\widehat{g} (n)|^2 =0. \label{simeqzero} \end{equation}
As we have seen above, (\ref{simeqzero}) is equivalent to
\begin{equation*}
\lim_{R\to \infty }\frac{1}{R} \sum_{n \in \stackrel{\circ}{D(R)} } |\widehat{f} (n)-\widehat{g} (n)|^2 =0 . 
\end{equation*}

Now let us consider the equation on ${\bf Z}^d $: \begin{equation}
(\widehat{H}-\lambda ) \widehat{u} =\widehat{f} . \label{section5_helmholtz} \end{equation}

%%%%%%%%%%%%%% Definition 5.5 %%%%%%%%%%%%%%%%%%

\begin{definition} \label{RadCond}
A solution $\widehat{u} (k)\in \widehat{\mathcal{B}}^{\ast} $ of (\ref{section5_helmholtz}) is said to be outgoing (for $+$) or incoming (for $-$) if it satisfies 
\begin{equation}
(\partial_{rad }\widehat{u} )(k) \simeq A_{\pm} (\lambda , \omega_k ) \widehat{u} (k) , \label{section5_radiation} 
\end{equation}
in the sense of (\ref{simeqzero}).
\end{definition}

%%%%%%%%%%%%%%%%%%%%%%% Theorem 5.6 %%%%%%%%%%%%%%%%%%%

\begin{theorem}
Let $\lambda \in I_d$.
If $\widehat{f} $ is compactly supported, $\widehat{R}(\lambda \pm i0 )\widehat{f} $ is an outgoing (for $+$) or incoming (for $-$) solution of the equation $(\widehat{H} -\lambda )\widehat{u} =\widehat{f}$.
\label{section5_radcondition} 
\end{theorem}
Proof. Since $x_{\infty } (\lambda , \omega_k )$ is homogeneous of degree $0$ in $k$ (see also the proof of Lemma 4.5), we have as $|k|\rightarrow \infty $ 
\begin{equation*}
x_{\infty }(\lambda , \omega_{k\pm e_j } )=x_{\infty }(\lambda , \omega_k ) +O(|k|^{-1} ) . 
\end{equation*}
Then we have for any fixed $n\in {\bf Z}^d $ 
\begin{equation}
e^{\pm in\cdot x_{\infty } (\lambda , \omega_{k\pm e_j} ) } -e^{\pm in \cdot x_{\infty } (\lambda , \omega_k )} =O(|k|^{-1} ) , \quad |k| \to \infty.
\label{lemma5_7_11}
\end{equation}
If $\widehat{f} $ is compactly supported, $\big( \widehat{\mathcal{G}}^{(\pm )} (\lambda ) \widehat{f} \big) (\omega_k ) $ is smooth with respect to $k$, so that we have from (\ref{lemma5_7_11})  \begin{equation}
\partial_{m-k} \left( a_{\pm} (\lambda , \omega_k) \big( \widehat{\mathcal{G}}^{(\pm )} (\lambda )\widehat{f} \big) (\pm \omega_k ) \right)=O(|k|^{-1} ). \label{lemma5_7_1} \end{equation}
We put $\widehat{u}^{(\pm)} =\widehat{R} (\lambda \pm i0 )\widehat{f} $.
 Theorem \ref{S4expansion_resol2} yields 
\begin{gather} 
\begin{split}
& (\partial_{rad} \widehat{u}^{(\pm )} )(k) \\
&=C_{\pm} |k|^{-(d-1 )/2} \left( \sum_{m\in \partial D(R(k)), m\sim k   } \big( \partial_{m-k} \Phi^{(\pm)}_{\lambda} \big) (k) \right) +O(|k|^{-(d+1)/2 } ) ,
\end{split} 
\label{lemma5_7_12}
\end{gather}
as $|k|\rightarrow \infty $,
where 
\begin{gather*}
C_{\pm} = \frac{1}{4} e^{\pm (3-d)\pi i /4}\sqrt{2\pi } ,\\
\Phi ^{(\pm)} _{\lambda} (k)= e^{\pm i k\cdot x_{\infty } (\lambda ,\omega_k ) }a_{\pm} (\lambda , \omega_k ) \big( \widehat{\mathcal{G}}^{(\pm)} (\lambda ) \widehat{f} \big) (\pm \omega_k ) . 
\end{gather*}
Lemma \ref{section5_estimate1} (2) and (\ref{lemma5_7_1}) imply the theorem. \qed

%%%%%%%%%%%%%%%%%%%%%%

\subsection{Rellich type theorem}

The following is an analogue of the Rellich type theorem for Schr{\"o}dinger operators in ${\bf R}^d$ (\cite{Re43}).

\begin{theorem} \label{Rellictype}
Let $\lambda \in (0,d) \setminus{\bf Z}$. Suppose a sequence $\{\widehat u(n)\}$ defined for $|n|\geq R_0 > 0$ satisfies
\begin{gather*}
(- \Delta_{disc} - \lambda)\widehat u = 0, \quad |n| > R_0, \\
\lim_{R\to\infty}\frac{1}{R}\sum_{R_0<|n|<R}|\widehat u(n)|^2 = 0.
\end{gather*}
Then there exists $R_1 > R_0$ such that $\widehat u(n) = 0$ for $|n| > R_1$.
\end{theorem}

For the proof, see \cite{IsMo}, Theorem 1.1.

\subsection{Uniqueness theorem}

%%%%%%%%%%%%%%%%%%%%%%%%%%% Theorem 5.7 %%%%%%%%%%%%%%%

\begin{theorem} 
Let $\lambda \in I_d$, and suppose that $\widehat{f} $ is compactly supported.
Let $\widehat{u}^{(\pm )} $ be the outgoing (for $+$) or incoming (for $-$) solution of the equation $(\widehat{H}-\lambda )\widehat{u}^{(\pm )} =\widehat{f} $. Then 
\begin{equation*}
(\widehat{u}^{(\pm )} , \widehat{f} )-(\widehat{f} , \widehat{u}^{(\pm )} )=2i \lim_{R \to \infty } \sum_{k\in \stackrel{\circ}{D(R)} \setminus \stackrel{\circ}{D(R-1)} } \mathrm{Im} A_{\pm } (\lambda , \omega_k)|\widehat{u}^{(\pm )} (k)|^2 . 
\end{equation*}
\label{section5_unique1}
 \end{theorem}

Proof. By Green's formula, we have 
\begin{gather*} 
\begin{split}
& \sum_{k \in \stackrel{\circ}{D(\rho )} } \left( (\Delta_{disc} \widehat{u} ^{(\pm )} )(k)\cdot \overline{\widehat{u}^{(\pm )} (k) } - \widehat{u}^{(\pm )} (k)\cdot \overline{(\Delta_{disc}\widehat{u}^{(\pm )} )(k) } \right) \\
&=\sum_{ k\in \partial D(\rho )} \left( (\partial _{\nu}^{D(\rho )} \widehat{u}^{(\pm )} )(k) \cdot \overline{\widehat{u} ^{(\pm )} (k)} - \widehat{u}^{(\pm )} (k) \cdot \overline{ (\partial_{\nu}^{D(\rho )} \widehat{u}^{(\pm )} )(k) } \right). 
\end{split} 
\end{gather*}
The left-hand side converges to $(\widehat{u}^{(\pm )} , \widehat{f} )-(\widehat{f}  , \widehat{u}^{(\pm )} )$ by the equation.
Changing the order of the summation,we can see that the right-hand side is equal to 
\begin{gather*} 
\begin{split}
&\frac{1}{4} \sum_{ k\in \partial D(\rho )} \sum_{m\in \stackrel{\circ}{D(\rho )}, m\sim k  } \left( \widehat{u} ^{(\pm)} (k) \cdot \overline{\widehat{u}^{(\pm )} (m)} - \widehat{u}^{(\pm)} (m) \cdot \overline{\widehat{u}^{(\pm )} (k)} \right) \\
&= \sum_{m\in \stackrel{\circ}{D(\rho )} \setminus \stackrel{\circ}{D(\rho -1)}} \left( (\partial _{rad}\widehat{u}^{(\pm )} )(m ) \cdot \overline{\widehat{u} ^{(\pm )} (m)} - \widehat{u} ^{(\pm )} (m) \cdot \overline{(\partial_{rad} \widehat{u}^{(\pm )} )(m)} \right) . 
\end{split} 
\end{gather*}
As $\rho \to \infty$, we can replace $\partial_{rad} \widehat{u}^{(\pm )}$ by $A_{\pm}(\lambda,\omega_k) \widehat{u}^{(\pm )}$, and prove the theorem. \qed

%%%%%%%%%% Theorem 5.8 %%%%%%%%%%%%%%%

\begin{theorem} \label{RadCondUnique}
Let $\lambda \in I_d$. 
If $\widehat{f} $ is compactly supported, then the outgoing solution of (\ref{section5_helmholtz}) is unique and given by $\widehat{R}(\lambda +i0 )\widehat{f} $.
The incoming solution is also unique and given by $\widehat{R} (\lambda -i0 )\widehat{f} $.
\label{section5_unique2}
\end{theorem}
Proof. In view of Theorem \ref{section5_radcondition}, we have only to prove the uniqueness.
Let $\widehat{u} $ be the outgoing solution of $(\widehat{H} -\lambda )\widehat{u} =0 $.
Then, by Theorem \ref{section5_unique1} and Lemma \ref{section5_estimate2} (2), we have $\lim_{R \to \infty } \sum_{k \in \stackrel{\circ}{D(R)} \setminus \stackrel{\circ}{D(R-1)}} |\widehat{u} (k)|^2 =0$.
This implies 
$$
\lim_{R \to \infty } \frac{1}{R} \sum_{k \in \stackrel{\circ}{D(R)}} |\widehat{u} (k)|^2 =0,
$$
i.e. $\widehat{u} \simeq 0$.
We can then use the Theorem \ref{Rellictype} and the unique continuation theorem (see \cite{IsMo}, Theorem 2.1) to see that $\widehat{u} =0$. \qed

%%%%%%%%%%%%%%%%%%%%%%%% Section 6 %%%%%%%%%%%%%%%%%%%%%%%

\section{Exterior problem}

%%%%%%%%%%%%%%%%%%%% subsection 6.1 %%%%%%%%%%%

\subsection{Helmholtz equation in an exterior domain}
Let $D(R)$ be a rectangular domain in (\ref{RectAngDom}), and take a sufficiently large integer $R_0 > 0$  such that
\begin{equation}
{\rm supp}\,\widehat V \subset \, \stackrel{\circ}{D(R_0)} .
\label{SuppVinD}
\end{equation}
We put
\begin{gather}
 \Omega_{int} = D(R_0 ),
\label{S6Omegaint} \\
\Omega_{ext} = {\bf Z}^d\setminus \stackrel{\circ}\Omega_{int}.
\label{S6Omegaext}
\end{gather}
Therefore $\stackrel{\circ}\Omega_{int} = \, \stackrel{\circ}{D(R_0)} \, = [-R_0 , R_0 ]^d \cap {\bf Z}^d $, and
\begin{equation}
\partial\Omega_{int} = \partial\Omega_{ext} = \bigcup_{j=1}^d \big\{
n\, ; \, |n_i|\leq R_0, \ (i\neq j), \ |n_j| = R_0+1 \big\}.
\label{S6partialOmegaintOmegaext}
\end{equation}

The spaces $\widehat{\mathcal B}$, $\widehat{\mathcal B}^{\ast}$ and $\widehat{\mathcal H}^s$ on $\stackrel{\circ}\Omega_{ext}$ are defined in the same way as in the whole space.  Let $\widehat H_{ext} = - \Delta_{disc}$ on $\Omega_{ext}$ with Dirichlet boundary condition.

%%%%%%%%%%%%%%%%%%%%%%%% Lemma 6.1 %%%%%%%%%%%%%%%%%%%%%

\begin{lemma} (1)
$\widehat H_{ext}$ is self-adjoint, and
$\sigma(\widehat H_{ext}) = [0,d].$ \\
\noindent
(2) $\ \sigma_p(\widehat H_{ext})\cap\left((0,d)\setminus{\bf Z}\right) = \emptyset.$
\label{section6_essspec}
\end{lemma}

Proof. The assertion (1) follows from the standard perturbation theory, and (2) is proved in Theorem 2.4 of \cite{IsMo}. \qed

\medskip
For the solution of the equation $(- \Delta_{disc}-\lambda)\widehat u = \widehat f$ in $\stackrel{\circ}\Omega_{ext}$, the radiation condition is defined in the same way as in \S 5. The following theorem is proved in the same way as in Theorem \ref{RadCondUnique}.  

%%%%%%%%%%%%%%%%%%%%%%%%%%%% Theorem 6.2 %%%%%%%%%%%%%%%%%%%%%%

\begin{theorem}
Let $\lambda \in I_d$. Then the solution of the equation $(- \Delta_{disc} - \lambda)\widehat u = 0$ in $\stackrel{\circ}\Omega_{ext}$, satisfying the Dirichlet boundary contidion and the outgoing (or incoming) radiation condition vanishes identically on $\Omega_{ext}$.
\label{S6extradiationcondition}
\end{theorem}

We prove the limiting absorption principle for
$\widehat R_{ext}(z) = (\widehat H_{ext}-z)^{-1}$.

%%%%%%%%%%%%%%%%%%%%%%%%%%% Theorem 6.3 %%%%%%%%%%%%%%%%%%%%%%%

\begin{theorem}
(1) For $\lambda \in I_d$ and $\widehat f \in \widehat{\mathcal B}$, the weak $\ast$-limit exists
$$
\lim_{\epsilon\to0}\widehat R_{ext}(\lambda \pm i\epsilon)\widehat f =:
\widehat R_{ext}(\lambda \pm i0)\widehat f \in \widehat{\mathcal B}^{\ast}.
$$
(2) For any compact set $J \subset I_d$, there exists a constant $C > 0$ such that
$$
\|\widehat R_{ext}(\lambda \pm i0)\widehat f\|_{\widehat{\mathcal B}^{\ast}} \leq C\|\widehat f\|_{\widehat{\mathcal B}}, \quad \lambda \in J.
$$
(3) For $\widehat f, \, \widehat g \in \widehat{\mathcal B}$, 
$$
I_d \ni \lambda \mapsto  \big( \widehat R_{ext}(\lambda \pm i0)\widehat f, \, \widehat g \big)
$$
is continuous. \\
(4) If $\widehat f$ is compactly supported, $\widehat R_{ext}(\lambda \pm i0)\widehat f$ satisfies the outgoing (for $+$) or incoming (for $-$) radiation condition.
 \end{theorem}

Proof. We prove the theorem for $\lambda + i0$. We extend $\widehat f \in \widehat{\mathcal{B}}$ and $\widehat u (z) = \widehat R_{ext}(z)\widehat f$ to be 0 outside $\Omega_{ext}$. Then it satisfies
$$
(\widehat H_0 - z)\widehat u(z) =  \widehat K\widehat u (z)  + \widehat f \quad \text{on} \quad {\bf Z}^d,
$$
where $\widehat K = \sum_n c_n\widehat P(n)$ is a finite sum of projections $\widehat P(n)$ to the site $n$.
Therefore
\begin{equation}
\widehat u(z) = \widehat R_0(z)\widehat K\widehat u(z) + \widehat R_0(z)\widehat f.
\label{S6resolventeq}
\end{equation}

Let $J$ be a compact set in $I_d$, and take $s > 1/2$. 
We first show that there exists a constant $C > 0$ such that
\begin{equation}
\|\widehat u(\lambda + i\epsilon ) \|_{\widehat{\mathcal H}^{-s}} \leq C\|\widehat f\|_{\widehat{\mathcal B}}, \quad \forall \lambda \in J, \quad \forall \epsilon > 0.
\label{S6Unifbound}
\end{equation}
 In fact, if this does not hold, there exists $z_{\mu } = \lambda_{\mu } + i\epsilon _{\mu }$, $\widehat f_{\mu } \in \widehat{\mathcal B}$, such that $\widehat u_{\mu } =\widehat R_{ext}(z_{\mu } ) \widehat f_{\mu } $ satisfies 
\begin{equation}
z_{\mu} \to \lambda \in J, \quad \|\widehat f_{\mu} \|_{\widehat{\mathcal B}}\to 0,  \quad \|\widehat u_{\mu } \|_{\widehat{\mathcal H}^{-s}} = 1 \quad \text{as} \quad \mu \to \infty .
\label{S6Absurd}
\end{equation}
One can then select a subsequence, which is denoted by $\{\widehat u_{\mu } \}$ again, such that $\widehat u_{\mu} $ converges weakly in $\widehat{\mathcal H}^{-s}$. Since $\widehat K$ is a finite dimensional operator, $\widehat K\widehat u_{\mu} $ converges in $\widehat{\mathcal B}$. Therefore, in view of (\ref{S6resolventeq}), we see that $\widehat u_{\mu }$ converges in $\widehat{\mathcal B}^{\ast}$, hence in $\widehat{\mathcal H}^{-s}$, to $\widehat u$ such that $\|\widehat u\|_{\widehat{\mathcal H}^{-s}}=1$. It satisfies
$$
(\widehat H_{ext}-\lambda)\widehat u = 0, \quad \widehat u = \widehat R_0(\lambda + i0)\widehat K\widehat u.
$$
Therefore $\widehat u$ is an outgoing solution. By Theorem \ref{S6extradiationcondition}, $\widehat u = 0$, which is a contradiction.

We next prove that for $s > 1/2$ and $\widehat f \in \widehat{\mathcal B}$, $\widehat R_{ext}(\lambda + i\epsilon)\widehat f$ converges strongly in $\widehat{\mathcal H}^{-s}$ as $\epsilon \to 0$. To prove it, we consider a sequence $u_{\mu} = \widehat R_{ext}(\lambda + i\epsilon_{\mu} )\widehat f$, $\epsilon_{\mu} \to 0$. 
Then by the same arguments as above, one can show that any subsequence of $\{u_{\mu} \}$ contains a sub-subsequence $\{u_{\mu '}\}$, which converges in $\mathcal H^{-s}$ to one and the same limit (independent of the choice of sub-subsequence). This proves the convergence of $\widehat R_{ext}(\lambda + i\epsilon)\widehat f$ as $\epsilon \to 0$.  Arguing similarly, one can also show that
$$
I_d \ni \lambda \mapsto  \widehat R_{ext}(\lambda + i0)\widehat f \in \widehat{\mathcal H}^{-s}
$$
is strongly continuous. The assertions of the theorem then follow from those for $\widehat R_0(\lambda + i0 )$ and the formula
$$
\widehat R_{ext}(\lambda + i0) = \widehat R_0(\lambda + i 0 ) \big( 1+ \widehat K
\widehat R_{ext}(\lambda + i0) \big).
\qed
$$

%%%%%%%%%%%%%%% subsection 6.2 %%%%%%%%%%%%%%%%%%%%%%

\subsection{Exterior and interior D-N maps}
Let $ \widehat{H}_{int}= -\Delta_{disc} +\widehat{V}$ be defined on $\Omega_{int}$ with Dirichlet boundary condition.
The {\it interior D-N map} is defined by 
\begin{equation}
\Lambda _{\widehat{V} } (\lambda ) \widehat{f} = \partial_{\nu}^{\Omega _{int}} \widehat{u}_{int}\Big|_{\partial\Omega_{int}},  \quad \lambda \not\in \sigma (\widehat{H}_{int}),
\label{IntDNmap}
\end{equation}
where $\widehat{u}_{int} $ is the solution of the equation 
\begin{equation}
(-\Delta _{disc} +\widehat{V} -\lambda )\widehat{u}_{int} =0 \quad \text{in} \quad \stackrel{\circ}\Omega_{int} , \quad \widehat{u} _{int}\Big|_{\partial \Omega _{int} } = \widehat{f}.
\label{S6uintequation}
\end{equation}

The {\it exterior D-N map} is defined by 
\begin{equation}
\Lambda _{ext}^{(\pm)} (\lambda ) \widehat{f} = - \partial _{\nu }^{\Omega_{ext} } \widehat{u}_{ext}^{(\pm)}\Big|_{\partial\Omega_{ext}},  \quad \lambda \in I_d,
\label{ExtDNmap}
\end{equation}
where $\widehat{u}_{ext}^{(\pm )} \in \widehat{\mathcal{B}}^{\ast } $ is the unique outgoing (for $+$) and incoming (for $-$) solution of the equation 
\begin{equation}
( \widehat{H}_{ext} - \lambda )\widehat{u}_{ext}^{(\pm)} =0 \quad \text{in} \quad \stackrel{\circ}\Omega _{ext} , \quad \widehat{u}_{ext}^{(\pm)} \Big|_{\partial \Omega _{ext} } = \widehat{f}. 
\label{S6uextpmequation}
\end{equation}
The existence of $\widehat{u} _{ext} ^{(\pm)} $ is shown by extending $\widehat{f}$ to be zero on  ${\bf Z}^d  \setminus \partial\Omega_{ext}$, and putting
\begin{equation*}
\widehat{u} _{ext} ^{(\pm)} = \widehat{f} - \widehat{R}_{ext} (\lambda \pm i0 ) ( -\Delta _{disc} - \lambda ) \widehat{f} . 
\end{equation*}
The uniqueness follows from Theorem 6.2.

\medskip
We  represent $\widehat u_{ext}^{(\pm)}$ in terms of exterior and interior D-N maps. In the following, for a subset $A$ in ${\bf Z}^d$, we use $\chi_{A}$ to mean either the operator of restriction
\begin{equation}
\chi_{A} : \ell^{\infty}({\bf Z}^d) \ni \widehat f \mapsto  \widehat f\Big|_{A},
\label{S6Restriction}
\end{equation}
or the operator of extension
\begin{equation}
\chi_{A} : \ell^{\infty}(A) \ni\widehat f \mapsto 
\left\{
\begin{array}{cl}
\widehat f, &  {\rm on} \quad A,\\
0 , & {\rm on} \quad {\bf Z}^d \setminus A,
\end{array}
\right.
\label{S6Extension}
\end{equation}
which will not confuse our argument. We  put
\begin{gather}
C(R_0) = \partial\Omega_{int} = \partial\Omega_{ext} = \Omega_{int}\cap\Omega_{ext},
\nonumber\\
\widehat{S}_{C(R_0)} = \frac{1}{4} \sum_{j=1}^d  \chi_{C(R_0)} \big( \widehat{S}_j + (\widehat{S}_j)^{\ast } \big) \chi_{C(R_0)},
\\
(\widehat{S}_j \widehat{u})(n) = \widehat{u} (n+ e_j ), \quad ( (\widehat{S}_j )^{\ast } \widehat{u} )( n)= \widehat{u} (n-e_j ),
\nonumber
\end{gather}
and also for $n \in C(R_0)$
\begin{equation}
\widetilde{\mathrm{deg}}_{C(R_0)}(n) = \frac{1}{4} \, ^{\#} \big\{ m \in C(R_0) \ ; \ |m-n |=1 \big\}.
\label{DefinceDeg}
\end{equation}
For $\lambda \in I_d\setminus\sigma(\widehat H_{int})$, we define the operator $ B_{C(R_0)}^{(\pm)} (\lambda ) \in {\bf B}(\ell^2(C(R_0)) )$ by 
\begin{equation}
 B_{C(R_0)}^{(\pm)} (\lambda ) = \Lambda _{\widehat{V}  }(\lambda ) - \Lambda _{ext}^{(\pm)}(\lambda )- \lambda +\frac{1}{4} \widetilde{\mathrm{deg}}_{C(R_0)}- \widehat{S}_{C(R_0)} ,
\label{DefineBOmega}
\end{equation}
where $\widetilde{\mathrm{deg}}_{C(R_0)}$ is the operator of multiplication by
$\widetilde{\mathrm{deg}}_{C(R_0)}(n)$.

%%%%%%%%%%%%%%% Lemma 6.4 %%%%%%%%%%%%%%%%%%%%%%%

\begin{lemma}
Assume that $\lambda \in I_d \setminus  \sigma(\widehat{H} _{int})$, $\widehat f \in \ell^2(C(R_0))$.
Let $\widehat u^{(\pm)}_{ext}$ and $\widehat{u}_{int}$ be the solutions of (\ref{S6uextpmequation}) and (\ref{S6uintequation}), respectively, and put 
$$
\widehat{u} ^{(\pm)} = \chi _{\stackrel{\circ}\Omega _{int}} \widehat{u}_{int} + \chi_{\stackrel{\circ}\Omega _{ext}} \widehat{u}_{ext}^{(\pm)} + \chi_{C(R_0)} \widehat{f}.
$$
Then we have
\begin{equation}
\widehat{u}^{(\pm)} (n) 
= ( \widehat{R} (\lambda \pm i0 )\chi_{C(R_0)} B_{C(R_0)}^{(\pm)} (\lambda ) \widehat{f} )(n), \quad n \in {\bf Z}^d. 
\label{S6upm(n)=RBf}
\end{equation}
In particular,
\begin{gather}
\widehat{u}_{ext}^{(\pm)} (n) 
= ( \widehat{R} (\lambda \pm i0 ) \chi_{C(R_0)} B_{C(R_0)}^{(\pm)} (\lambda ) \widehat{f} )(n), \quad n \in \, \stackrel{\circ}\Omega_{ext}, 
\label{S6upmext(n)=RBf} \\
\widehat{f} (n) 
= ( \widehat{R} (\lambda \pm i0 ) \chi_{C(R_0)} B_{C(R_0)}^{(\pm)} (\lambda ) \widehat{f} ) (n), \quad n \in C(R_0). 
\label{S6hatf=RBf}
\end{gather}
\label{exteriorsol} 
\end{lemma}
Proof. Let $\widehat{r}(n,m;\lambda \pm i0 )$ be the resolvent kernel, i.e. 
$$
\widehat{r}(n,m;\lambda \pm i0 )=\left(\widehat{R} (\lambda \pm i0)\widehat \delta_{m}\right)(n),
$$
 where $\widehat \delta_m(n) = \delta_{mn}$.
As in the proof of Theorem \ref{section5_unique1}, by Green's formula,
\begin{gather} 
\begin{split}
&\sum_{ n \in (\stackrel{\circ}\Omega_{int} \cup \stackrel{\circ}\Omega_{ext})\cap \stackrel{\circ}{D(R)} } \left( (\Delta_{disc} \widehat{u}^{(\pm)} )(n) \widehat{r} (n,m;\lambda \pm i0)- \widehat{u}^{(\pm)} (n) ( \Delta_{disc} \widehat{r})(n,m;\lambda \pm i0 ) \right) \\
&= \sum_{ n \in \partial \Omega _{int} } \left( ( \partial _{\nu} ^{\Omega_{int} } \widehat{u}^{(\pm)} )(n) \widehat{r} (n,m;\lambda \pm i0 ) - \widehat{u}^{(\pm)} (n) ( \partial _{\nu}^{\Omega_{int}}\widehat{r})(n,m;\lambda \pm i0) \right) \\
& \quad + \sum_{ n\in \partial \Omega _{ext} } \left( (\partial _{\nu}^{\Omega_{ext} } \widehat{u}^{(\pm)} )(n) \widehat{r} (n,m;\lambda \pm i0 ) - \widehat{u}^{(\pm)} (n) (\partial _{\nu}^{\Omega_{ext} } \widehat{r} )(n,m;\lambda \pm i0 ) \right) \\
&\quad + \sum_{n\in \stackrel{\circ}{D(R)} \setminus \stackrel{\circ}{D(R-1)} } \left( (\partial _{rad} \widehat{u}^{(\pm )} )(n) \widehat{r} (n,m;\lambda \pm i0 ) - \widehat{u}^{(\pm)} (n) (\partial _{rad} \widehat{r})(n,m;\lambda \pm i0 ) \right),
 \label{greencal} 
\end{split} 
\end{gather}
for sufficiently large integer $R>0$. 
By the equations (\ref{S6uintequation}) and (\ref{S6uextpmequation}), the left-hand side of (\ref{greencal}) is equal to 
\begin{gather} \begin{split}
&\sum_{ n \in (\stackrel{\circ}\Omega _{int} \cup \stackrel{\circ}\Omega_{ext} )\cap \stackrel{\circ}{D(R)}} \widehat{u}^{(\pm)} (n) (( -\Delta_{disc} +\widehat{V} -\lambda ) \widehat{r} )(n,m;\lambda \pm i0 ) \\
&=   \sum_{ n \in ( \stackrel{\circ}\Omega _{int} \cup \stackrel{\circ}\Omega_{ext} )\cap \stackrel{\circ}{D(R)} } \widehat{u}^{(\pm)}(n) \delta_{nm},
\label{S6LHSGreen} \end{split}
\end{gather}
for any $m \in {\bf Z}^d$.
Note that, by our definitions of $ \Lambda _{\widehat{V}} (\lambda ) $ and $\Lambda _{ext} ^{(\pm)} (\lambda )$, 
\begin{gather*} 
 \partial _{\nu }^{\Omega_{int} } \widehat{u} ^{(\pm)} = \partial _{\nu } ^{\Omega_{int} } \widehat{u}_{int} = \Lambda _{\widehat{V}} (\lambda ) \widehat{f}, \\
 \partial _{\nu} ^{\Omega_{ext} } \widehat{u} ^{(\pm)} = \partial_{\nu} ^{\Omega_{ext} } \widehat{u}_{ext}^{(\pm)} = - \Lambda _{ext}^{(\pm)} (\lambda )\widehat{f} . 
\end{gather*}
The sum $\sum_{n \in \partial \Omega_{int} } +\sum_{n\in \partial \Omega_{ext} } $ in the right-hand side of (\ref{greencal}) is then equal to 
\begin{gather} 
\begin{split}
&\sum_{ n \in C(R_0)} \left( ( \Lambda _{\widehat{V}  } (\lambda ) \widehat{f} ) (n) \widehat{r } (n,m;\lambda \pm i0 )- \widehat{f} (n) (\partial _{\nu }^{\Omega_{int} } \widehat{r} )(n,m;\lambda \pm i0 ) \right) \\
&\quad - \sum_{ n\in C(R_0) } \left(( \Lambda _{ext}^{(\pm )} (\lambda ) \widehat{f} )(n) \widehat{r } (n,m;\lambda \pm i0 ) + \widehat{f} (n) ( \partial _{\nu}^{\Omega_{ext} } \widehat{r })(n,m;\lambda \pm i0) \right) \\
&=\sum_{ n \in C(R_0)} \widehat{r} ( n,m; \lambda \pm i0 )\chi _{C(R_0)}(n)\left((\Lambda _{ \widehat{V} }(\lambda ) -  \Lambda _{ext}^{(\pm )}  (\lambda ) ) \widehat{f}\right) (n) \\
& \quad - \sum _{ n \in C(R_0)}  \widehat{f} (n) \left( ( \partial _{\nu} ^{\Omega _{int} } + \partial _{\nu} ^{\Omega _{ext} } ) \widehat{r} \right)( n,m ; \lambda \pm i0 ). 
\end{split} 
\label{S6Longformula}
\end{gather}
For $n \in C(R_0)$,
\begin{equation}
\begin{split}
& \left((\partial_{\nu}^{\Omega_{int}} + \partial_{\nu}^{\Omega_{ext}})\widehat r \right)(n,m;\lambda \pm i0) \\
= & \frac{1}{4}\sum_{k \in \stackrel{\circ}\Omega_{int}\cup\stackrel{\circ}\Omega_{ext} , k \sim n }
\big(\widehat r(n,m;\lambda \pm i0) - \widehat r (k,m;\lambda \pm i0)\big) \\
=& - (\Delta_{disc}\widehat r)(n,m;\lambda \pm i0) - 
\frac{1}{4}\sum_{k \in C(R_0), k\sim n}\big(\widehat r(n,m;\lambda \pm i0) - \widehat r (k,m;\lambda \pm i0)\big).
\end{split}
\nonumber
\end{equation}
Therefore, the second term of the right-hand side of (\ref{S6Longformula}) is computed as follows:
\begin{gather*} 
\begin{split} 
&-  \sum_{ n \in C(R_0)} \widehat{f} (n) \ (-\Delta_{disc} \widehat{r}) (n,m ; \lambda \pm i0 )  \\
& \quad +  
\frac{1}{4} \sum_{ n \in C(R_0) } \widehat{f} (n) \Big( \sum_{k\in 
C(R_0), k \sim n } \big(\widehat{r} (n,m;\lambda \pm i0 )- \widehat{r} ( k , m ; \lambda \pm i0 ) \big) \Big) \\
&=-\sum_{n\in C(R_0)} \widehat{f}(n) \delta_{nm} +  \sum_{ n \in C(R_0) } \big( -\lambda + \frac{1}{4} \widetilde{\mathrm{deg} }_{C(R_0)} (n) \big) \widehat{f} (n) \ \widehat{r} (n,m;\lambda \pm i0 ) \\
& \quad - \frac{1}{4} \sum _{k \in C(R_0 )} \widehat{r} (k,m;\lambda \pm i0 ) \sum_{n \in C(R_0), n \sim k} \widehat{f} (n), 
 \end{split} 
\end{gather*}
where, in the 3rd line, we have  used the fact that 
\begin{equation*}
\big( (-\Delta_{disc} -\lambda ) \widehat{r} \big) (n,m;\lambda \pm i0 ) =\delta_{nm}  , \quad n, m\in {\bf Z}^d,
 \end{equation*}
and exchanged the order of summation in the 4th line.
Note that
$$
\sum_{n \in C(R_0), n\sim k}\widehat f(n) = \sum_{j=1}^d\big((\widehat S_j + (\widehat S_j)^{\ast})\chi_{C(R_0)}\widehat f\big)(k).
$$
Since we have for any $m \in \, \stackrel{\circ}{D(R)}$ 
\begin{equation*}
\sum_{n\in (\stackrel{\circ}\Omega_{int} \cup \stackrel{\circ}\Omega_{ext} )\cap \stackrel{\circ}{D(R)} } \widehat{u}^{(\pm )} (n)\delta_{nm} +\sum_{n\in C(R_0)} \widehat{f}(n)\delta_{nm } =\widehat{u}^{(\pm )} (m), 
\end{equation*}
(\ref{greencal}) turns out to be
\begin{gather*} 
\begin{split}
 &\widehat{u}^{(\pm )}(m)  
=\sum_{ n \in C(R_0) } \widehat{r} ( n,m; \lambda \pm i0 )\big((\Lambda _{ \widehat{V} }(\lambda ) - \Lambda _{ext}^{(\pm )}  (\lambda ) ) \widehat{f}\big) (n) \\
& \ \ \ \ \ \ \ \ \ \ \ \ \ \ +  \sum_{ n \in C(R_0) } \big(- \lambda + \frac{1}{4} \widetilde{\mathrm{deg} }_{C(R_0)} (n) \big) \widehat{f} (n) \ \widehat{r} (n,m;\lambda \pm i0 ) \\
&\ \ \ \ \ \ \ \ \ \ \ \ \ \ -\frac{1}{4} \sum_{n\in C(R_0 )} \widehat{r} (n,m;\lambda \pm i0 ) \sum_{j=1}^d\big((\widehat S_j + (\widehat S_j)^{\ast})\chi_{C(R_0)}\widehat f\big) (n) \\
&\quad + \sum_{n\in \stackrel{\circ}{D(R)} \setminus \stackrel{\circ}{D(R-1)}}\left( (\partial _{rad} \widehat{u}^{(\pm )} )(n) \widehat{r}(n,m;\lambda \pm i0 ) - \widehat{u}^{(\pm)}(n) (\partial _{rad} \widehat{r})(n,m;\lambda \pm i0 ) \right)  , 
\end{split} 
\end{gather*}
for any $m\in \, \stackrel{\circ}{D(R) }$.
In view of (\ref{DefineBOmega}), we have thus arrived at
\begin{gather*} 
\begin{split}
&\widehat{u}^{(\pm )} (m)  
=  \big( \widehat{R}(\lambda \pm i0 ) \chi_{C(R_0)} B^{(\pm )} _{C(R_0 )} (\lambda ) \widehat{f} \big) (m) \\
&\quad + \sum_{n\in  \stackrel{\circ}{D(R)} \setminus \stackrel{\circ}{D(R-1)}}\left( (\partial _{rad} \widehat{u}^{(\pm )} )(n) \widehat{r} (n,m;\lambda \pm i0 ) - \widehat{u}^{(\pm)}(n)(\partial _{rad} \widehat{r} )(n,m;\lambda \pm i0 ) \right)  . 
\end{split} 
\end{gather*} 
Taking the average of the sum with respect to $R$ in the above equality, we have 
\begin{gather} 
\begin{split} \label{S6pre}
\widehat{u}^{(\pm )}(m) & =  \big( \widehat{R}(\lambda \pm i0 ) \chi_{C(R_0)} B^{(\pm )} _{C(R_0 )} (\lambda ) \widehat{f} \big) (m) \\
& + \frac{1}{R} \sum_{n\in \stackrel{\circ}{D(R)}} \big( (\partial _{rad} -A_{\pm} (\lambda , \omega_n ) ) \widehat{u}^{(\pm )} \big) (n) \widehat{r} (n,m;\lambda  \pm i0 ) \\
& -\frac{1}{R} \sum_{n \in \stackrel{\circ}{D(R)}} \widehat{u}^{(\pm )} (n) \big( ( \partial _{rad} -A_{\pm }(\lambda , \omega_n )) \widehat{r} \big)(n,m;\lambda \pm i0 ) , 
\end{split} 
\end{gather}
up to a term of $O(R^{-1})$.
By the radiation condition, we have
 \begin{gather*} \begin{split} 
\frac{1}{R} & \left| \sum_{n\in \stackrel{\circ}{D(R)}}  \big( (\partial _{rad} -A_{\pm} (\lambda , \omega_n ) ) \widehat{u}^{(\pm )} \big) (n) \widehat{r} (n,m;\lambda  \pm i0 ) \right| \\
& \leq \left( \frac{1}{R} \sum_{n \in \stackrel{\circ}{D(R)}} \big| \big( ( \partial _{rad} -A_{\pm } (\lambda ,\omega_n) ) \widehat{u}^{(\pm )} \big) (n) \big|^2 \right)^{1/2} \\
& \quad \times\left( \frac{1}{R} \sum_{n \in \stackrel{\circ}{D(R)}} | \widehat{r} (n,m;\lambda \pm i0 )|^2 \right)^{1/2},
 \end{split} \end{gather*}
which tends to zero as $R \rightarrow \infty$.
The third term of the right-hand side of (\ref{S6pre}) is estimated similarly.
This proves the lemma.
\qed 

\medskip

%%%%%%%%%%%%%%%%% Lemma 6.5 %%%%%%%%%%%%%%%%%%%%%%

\begin{lemma}
Suppose $\lambda \in I_d \setminus  \sigma (\widehat{H}_{int})$.
Then for $\widehat{f} , \, \widehat{g} \in \ell^2 (C(R_0) )$, we have 
\begin{gather}
\big( \Lambda _{\widehat{V} } (\lambda ) \widehat{f} , \, \widehat{g} \big)_{ \ell^2 (C(R_0) ) }
 = \big(  \widehat{f} , \, \Lambda _{\widehat{V} } (\lambda ) \widehat{g} \big)_{ \ell^2 (C(R_0) ) }, 
\label{adjointdnint} \\
\big( \Lambda _{ext } ^{(\pm)} (\lambda ) \widehat{f} , \, \widehat{g} \big)_{ \ell^2 (C(R_0) ) }
 = \big(  \widehat{f} , \, \Lambda _{ext}^{(\mp)} (\lambda ) \widehat{g} \big)_{ \ell^2 (C(R_0) ) } . \label{adjointdnext} 
\end{gather} 
\label{adjointlem}
\end{lemma}
Proof. The first equality (\ref{adjointdnint}) follows from Green's formula.
We shall prove (\ref{adjointdnext}).
Let $\widehat{u} $ be the outgoing solution of (\ref{S6uextpmequation}), 
and $\widehat{v} $  the incoming solution of (\ref{S6uextpmequation}) with $\widehat f$ replaced by $\widehat g$.
For a sufficiently large integer $R>0$, we have by Green's formula
\begin{gather*} 
\begin{split}
0&= \sum_{ n \in (D(R)\cap\Omega_{ext})^{\circ}} \left( ((\widehat{H}_{ext} -\lambda )\widehat{u})(n) \cdot \overline{\widehat{v} (n)} - \widehat{u}(n) \cdot \overline{ ( (\widehat{H}_{ext}-\lambda )\widehat{v})  (n)} \right) \\
&= \sum_{n\in \partial D(R)} \left( -( \partial _{\nu}^{D(R)} \widehat{u} )(n) \cdot \overline{\widehat{v} (n)} + \widehat{u}(n) \cdot \overline{ (\partial_{\nu}^{D(R)} \widehat{v})(n)} \right) \\
&\quad +\sum_{n\in \partial \Omega_{ext} }  \left( -( \partial _{\nu}^{\Omega_{ext} } \widehat{u} )(n) \cdot \overline{\widehat{v} (n)} + \widehat{u}(n) \cdot \overline{ (\partial_{\nu}^{\Omega_{ext}} \widehat{v})(n)} \right) . 
\end{split} \end{gather*}
As in the proof of Theorem \ref{section5_unique1}, we have
 \begin{gather*} 
\begin{split}
&\sum_{ n \in \partial D(R)} \left( -( \partial _{\nu}^{D(R)} \widehat{u} )(n) \cdot \overline{\widehat{v} (n)} + \widehat{u}(n) \cdot \overline{ (\partial_{\nu}^{D(R)} \widehat{v})(n)} \right) \\
&= \sum_{n\in \stackrel{\circ}{D(R)} \setminus \stackrel{\circ}{D(R-1)}} \left( -(\partial _{rad}\, \widehat{u} )(n) \cdot \overline{ \widehat{v} (n)} + \widehat{u} (n) \cdot \overline{(\partial_{rad} \, \widehat{v} (n) } \right) . 
\end{split} 
\end{gather*}
This implies
\begin{gather*} 
\begin{split}
& \big( \Lambda_{ext}^{(+)} (\lambda )\widehat{f} , \, \widehat{g} \big)_{\ell^2 (\partial \Omega_{ext} ) } - \big(\widehat{f} , \, \Lambda _{ext}^{(-)} (\lambda )\widehat{g} \big)_{\ell^2 (\partial \Omega_{ext} )} \\
&= \sum_{ n \in \stackrel{\circ}{D(R)} \setminus \stackrel{\circ}{D(R-1)}} \big( ( \partial _{rad}\, \widehat{u})(n) - A_+ (\lambda , \omega_n )\widehat{u} (n) \big)  \overline{\widehat{v} (n)} \\
& \quad - \sum_{ n \in \stackrel{\circ}{D(R)} \setminus \stackrel{\circ}{D(R-1)}}  \widehat{u} (n)\overline{\big( ( \partial _{rad}\, \widehat{v})(n) - A_- (\lambda , \omega_n )\widehat{v} (n) \big)}. 
\end{split} 
\end{gather*}
Then, taking the average of the sum with respect to $R$, we have
 \begin{gather*}
 \begin{split}
& \big( \Lambda_{ext}^{(+)} (\lambda )\widehat{f} , \, \widehat{g} \big)_{\ell^2 (\partial \Omega_{ext} ) } - \big( \widehat{f} , \, \Lambda _{ext}^{(-)} (\lambda )\widehat{g} \big)_{\ell^2 (\partial \Omega_{ext} )} \\
&= \frac{1}{R} \sum_{n\in \stackrel{\circ}{D(R)}} \big( ( \partial _{rad}- A_+ (\lambda , \omega_n ))\widehat{u} (n) \big)  \overline{\widehat{v} (n)}  \\
&
- \frac{1}{R} \sum_{n\in \stackrel{\circ}{D(R)}}  \widehat{u} (n) \overline{\big( ( \partial _{rad} - A_- (\lambda , \omega_n ))\widehat{v} (n)\big) },
 \end{split} 
\end{gather*}
up to a term of $O(R^{-1})$. 
By the radiation condition, we can see that the right-hand side tends to zero as $R\rightarrow \infty $ as in the estimate of (\ref{S6pre}).
This proves (\ref{adjointdnext}). \qed

%%%%%%%%%%%%%%%%%%%%%%%%%%%%% Section 7 %%%%%%%%%%%%%%%%%%%%

\section{Scattering amplitude and D-N maps}

\subsection{Far-field pattern}
We introduce the operator $ \widehat{\Gamma }^{(\pm)} (\lambda ) $  by 
\begin{equation}
 \widehat{ \Gamma }^{(\pm)} (\lambda )= \widehat{\mathcal{G}} ^{(\pm)} (\lambda ) \chi_{C(R_0)} B_{C(R_0) }^{(\pm)}(\lambda ) : \ell^2 (C(R_0) )\rightarrow L^2 (S^{d-1 }) .
\label{S6Gammapmdef}
\end{equation}
The main purpose of this subsection is to show that $ \widehat{\Gamma }^{(\pm)} (\lambda ) $ is 1 to 1 (Lemma \ref{neartofarprop}). 

Although defined through $\widehat{\mathcal G}^{(\pm)}(\lambda)$, $ \widehat{\Gamma }^{(\pm)} (\lambda ) $ does not depend on $\widehat V$. It is seen by the  next lemma which follows from Lemma \ref{exteriorsol} and Theorem \ref{Rlambdapmi0Asympto}.

%%%%%%%%%%%%%% Lemma 7.1 %%%%%%%%%%%%%%%%%%%%%%%%%%

\begin{lemma}
Suppose $\lambda \in I_d \setminus  \sigma (\widehat{H}_{int})$.
Let $\widehat{u}_{ext}^{(\pm)} $ be the solution of (\ref{S6uextpmequation}).
Then we have 
\begin{gather*} 
\begin{split}
\widehat{u}_{ext}^{(\pm )}(k) 
&= e^{ \pm (3-d) \pi i /4 } \sqrt{2 \pi } |k| ^{ -(d-1)/2} e^{ \pm ik \cdot x_{\infty} (\lambda , \omega _k ) } a_{ \pm } (\lambda , \omega _k ) \big( \widehat{\Gamma } ^{(\pm )} (\lambda ) \widehat{f} \big) (\pm \omega _k ) \\
& \quad + O(|k|^{- (d+1)/2} )
\end{split} 
\end{gather*}
as $|k| \rightarrow \infty $.
\label{exteriorasymp}
 \end{lemma}

We need resolvent equations for $\widehat{R}_{ext} (\lambda \pm i0 )$. Note that by (\ref{DefineBOmega}) and Lemma \ref{adjointlem}
\begin{equation*} 
( B_{\partial \Omega }^{(\pm)} (\lambda ) )^{\ast} = \Lambda _{\widehat{V}}(\lambda ) -\Lambda _{ext}^{(\mp )}(\lambda )- \lambda +\frac{1}{4} \widetilde{\mathrm{deg} }_{C(R_0)} -\widehat{S}_{C(R_0)} = B_{\partial \Omega }^{(\mp)} (\lambda ).
\end{equation*}

%%%%%%%%%%%%%%% Lemma 7.2 %%%%%%%%%%%%%%%%%%

\begin{lemma} \label{Rextkinfty}
\begin{gather} 
\widehat{R}_{ext} (\lambda \pm i0 )= \widehat{R} _0 (\lambda \pm i0 )- \widehat{R} (\lambda \pm i0 ) \chi_{C(R_0)} B_{\partial \Omega }^{(\pm)}(\lambda ) \chi_{C(R_0)}\widehat{R}_0 (\lambda \pm i0 ) . \label{reseq} \\
\widehat{R}_{ext} (\lambda \pm i0 )= \widehat{R} _0 (\lambda \pm i0 )- \widehat{R} _0 (\lambda \pm i0 )\chi_{C(R_0)}  B_{\partial \Omega }^{(\pm)} (\lambda ) \chi_{C(R_0)} \widehat{R} (\lambda \pm i0 ). \label{reseqadj} 
\end{gather}
\end{lemma}
Proof. Since $\widehat{v}_0 = \widehat{R} (\lambda \pm i0 ) \chi_{C(R_0)} B_{C(R_0)}^{(\pm )}(\lambda ) \chi_{C(R_0)}\widehat{R}_0 (\lambda \pm i0 ) \widehat{f}$ satisfies the equation 
\begin{equation*}
(\widehat{H}_{ext} - \lambda )\widehat{v}_0 =0 \quad \text{in} \quad \stackrel{\circ}\Omega_{ext} , \quad \widehat{v}_0 |_{ \partial \Omega_{ext} } = \widehat{R}_0 (\lambda \pm i0 ) \widehat{f} , 
\end{equation*}
 we have (\ref{reseq}). 
Taking the adjoint, we obtain (\ref{reseqadj}). \qed

\medskip
We introduce the generalized Fourier transform in the exterior domain. 
We put 
\begin{equation*}
\widehat{ \mathcal{F} }_{ext}^{(\pm )} (\lambda ) = \widehat{\mathcal{F} }_0 (\lambda )\left( 1- \chi_{C(R_0)}B_{C(R_0) } ^{(\pm)} (\lambda) \chi_{C(R_0)}  \widehat{R} (\lambda \pm i0)\right) ,\quad \lambda \in I_d \setminus \sigma_p (\widehat{H}_{int}),
\end{equation*}
and, in the same way as (\ref{S4G0pmjlambdadefine}), we define 
$$
( \widehat{ \mathcal{G}} ^{(\pm)}_{ext} (\lambda ) \widehat{f} )( \omega )= ( \widehat{\mathcal{F}}_{ext} ^{(\pm)} (\lambda )\widehat{f} )( \theta (\lambda , \omega )).
$$
Lemmas \ref{R0lambdakinfty} and \ref{Rextkinfty} imply that  as $|k| \rightarrow \infty $, 
\begin{gather*} 
\begin{split} 
& ( \widehat{R}_{ext} (\lambda \pm i0 ) \widehat{f}) (k) \\ 
&= e^{ \pm (3-d) \pi i / 4} \sqrt{ 2 \pi } |k| ^ { -(d-1)/2 } e^{ \pm i k \cdot x_{\infty } (\lambda , \omega _k )} a_{\pm} (\lambda , \omega _k ) ( \widehat{\mathcal{G}} ^{(\pm )}_{ext} (\lambda ) \widehat{f} )(\pm \omega _k ) \\
 & \ \ \ \ \ \ \ \ \ \ +O(|k| ^{-(d+1 )/2} ). 
\end{split} 
\end{gather*} 
This formula shows that $\widehat{ \mathcal{G}} ^{(\pm)}_{ext} (\lambda )$ does not depend on $\widehat V$.

%%%%%%%%%%%%%%%%%%%%%%%%%% Lemma 7.3 %%%%%%%%%%%%%%%%%%%%%%%%55

\begin{lemma}
For any $ \widetilde{\phi} \in L^2 ( S^{d-1} )$, $\widehat{ \mathcal{G} }_{ext}^{(- )} (\lambda )^{\ast } \widetilde{\phi } $ satisfies the equation 
\begin{equation*}
( \widehat{H}_{ext} - \lambda ) \widehat{ \mathcal{G} }_{ext}^{(- )} (\lambda )^{\ast } \widetilde{\phi }=0 \quad \text{in} \quad \stackrel{\circ}\Omega_{ext} , \quad \Big(\widehat{ \mathcal{G} }_{ext}^{(- )} (\lambda )^{\ast } \widetilde{\phi }\Big)\Big|_{ \partial \Omega_{ext} } =0 , 
\end{equation*}
and $\widehat{ \mathcal{G} }_{ext}^{(-)} (\lambda )^{\ast } \widetilde{\phi } - \widehat{ \mathcal{G} }_0 (\lambda )^{\ast } \widetilde{\phi }$ is outgoing.
\label{S7exteriorfourier}
\end{lemma}

Proof. By the definition, we have 
\begin{equation*}
\widehat{ \mathcal{G} }_{ext}^{(-)} (\lambda )^{\ast } \widetilde{\phi } = \left( 1 - \widehat{R}(\lambda +i0 ) \chi_{C(R_0)} B_{C(R_0) }^{(+)}(\lambda ) \chi_{C(R_0)}\right) \widehat{ \mathcal{G} }_0 (\lambda )^{\ast } \widetilde{\phi } . 
\end{equation*}
By Lemma \ref{exteriorsol},
$ \widehat{v}=\widehat{R}(\lambda +i0 )\chi_{C(R_0)}  B_{C(R_0) }^{(+)} (\lambda ) \chi_{C(R_0)}\widehat{ \mathcal{G} }_0 (\lambda )^{\ast } \widetilde{\phi } $ satisfies the equation 
\begin{equation*}
( \widehat{H}_{ext} - \lambda )\widehat{v} =0 \quad \text{in} \quad \stackrel{\circ}\Omega_{ext} , \quad \widehat{v} |_{ \partial \Omega_{ext}} = \widehat{\mathcal{G}}_0 (\lambda )^{\ast } \widetilde{\phi } . 
\end{equation*}
The lemma then follows if we note that $\widehat{\mathcal{G}}_0 (\lambda )^{\ast } \widetilde{\phi } $ satisfies 
\begin{equation*}
( \widehat{H}_{ext} - \lambda ) \widehat{\mathcal{G}}_0 (\lambda )^{\ast } \widetilde{\phi } =0 \quad \text{in} \quad \stackrel{\circ}\Omega_{ext} .
\qed
 \end{equation*}

%%%%%%%%%%%%%%%%%%%%%%%%% Lemma 7.4 %%%%%%%%%%%%%%%%%%%

\begin{lemma}
Suppose $\lambda \in I_d \setminus \sigma (\widehat{H}_{int} )$. \\
 (1) \ $\widehat{\Gamma}^{(\pm)} (\lambda ) : \ell^2(C(R_0)) \to L^2(S^{d-1})$ is 1 to 1. \\
 \noindent
 (2) \ $\widehat{\Gamma}^{(\pm)} (\lambda )^{\ast} : L^2(S^{d-1}) \to \ell^2(C(R_0)) $ is onto.
\label{neartofarprop}
\end{lemma}
Proof. Let us show (1).
Suppose $\widehat{\Gamma}^{(\pm )} (\lambda ) \widehat{f} =0 $ and let $\widehat{u}_{ext}^{(\pm )} $ be the solution of (\ref{S6uextpmequation}).
From Lemma \ref{exteriorasymp} and the assumption, we have $\widehat{u}_{ext}^{(\pm )} \simeq 0$.
Then we see that $\widehat{u}_{ext}^{(\pm )} $ is compactly supported by Theorem \ref{Rellictype}. By the unique continuation property (see \cite{IsMo}, Theorem 2.3),  we then obtain $\widehat{f} =0$, which proves (1). This implies that the range of $\widehat{\Gamma}^{(\pm )} (\lambda )^{\ast}$ is dense. Since $\ell^2(C(R_0))$ is finite dimensional, (2) follows. \qed

%%%%%%%%%%%%%%%%%% subsection 7.2 %%%%%%%%%%%%%

\subsection{Scattering amplitude}
Recall that the scattering amplitude in the whole space is defined by (\ref{S3ScattAmpWholeSp}). Passing to $M_{\lambda}$, we rewrite it as
\begin{equation}
\mathcal{A} (\lambda )= \widehat{\mathcal{G}} ^{(+)} (\lambda ) \widehat{V} \widehat{ \mathcal{G}}_0 (\lambda )^{\ast }.
\label{S7ScattAmpWholeReparamet}
\end{equation}
The scattering amplitude for the exterior domain is defined by 
\begin{equation}
A_{ext} (\lambda )= \widehat{ \mathcal{F}} ^{(+)} (\lambda ) \chi_{C(R_0)} B_{C(R_0) } ^{ (+)} (\lambda ) \chi_{C(R_0)}  \widehat{ \mathcal{F}}_0 (\lambda ) ^{\ast }.
\label{S7ScattAmpExt} 
\end{equation}
As in the case of ${\bf Z}^d$, we use its reparametrization on $M_{\lambda}$:
\begin{equation}
\mathcal{A} _{ext} (\lambda )= \widehat{ \mathcal{G}} ^{(+)} (\lambda )\chi_{C(R_0)} B_{C(R_0 )} ^{ (+)} (\lambda ) \chi_{C(R_0)} \widehat{ \mathcal{G}}_0 (\lambda ) ^{\ast }.
\label{S7ScattAmpExtReparamet}
\end{equation}
Then we have as $|k|\rightarrow \infty $ 
\begin{gather} 
\begin{split} 
& ( \widehat{ \mathcal{G} }_{ext}^{(-)} (\lambda )^{\ast } \widetilde{\phi })(k) - (\widehat{ \mathcal{G} }_0 (\lambda )^{\ast } \widetilde{\phi } )(k) \\
&= - e^{ (3-d ) \pi i /4} \sqrt{2 \pi } |k| ^{ - (d-1 )/2 } e^{ ik \cdot x _{\infty} (\lambda , \omega _k )} a_{+ } (\lambda , \omega _k ) ( \mathcal{A}_{ext} (\lambda )\widetilde{\phi } )(\omega _k ) \\
&\quad +O(|k| ^ {-(d+1)/2} ) .
\label{G-lambda-G0labdaAsympto}
\end{split} 
\end{gather}
In fact, the left-hand side  is equal to
$ - \widehat R(\lambda + i0) \chi_{C(R_0)} B^{(+)}_{C(R_0)}(\lambda) \chi_{C(R_0)}\widehat{\mathcal  G}_0(\lambda)^{\ast} \widetilde{\phi }$. Using  Theorem \ref{Rlambdapmi0Asympto}, we obtain (\ref{G-lambda-G0labdaAsympto}).
 
%%%%%%%%%%%% subsection 7.3 %%%%%%%%%%%%%%%%

\subsection{Single layer and double layer potentials} 
We have already introduced the operator $\widehat R(\lambda \pm i0)\chi_{C(R_0)} B^{(\pm)}_{C(R_0)}(\lambda)$, which is an analogue of the double layer potential. We also need a counter part for the single layer potential, which is an operator on $\ell^2 (C(R_0) )$ defined by
\begin{equation*}
M_{C(R_0)}^{(\pm)} (\lambda )\widehat{f} = \left(\widehat{R} (\lambda \pm i0 )
\chi_{C(R_0) }\widehat{f}\right)\Big |_{C(R_0)}  
\end{equation*}
for $\widehat{f} \in \ell^2 (C(R_0) )$. 

The following lemma is a direct consequence of (\ref{S6hatf=RBf}) and the fact that $M_{C(R_0)}^{(\pm )} (\lambda )$ corresponds to $\chi_{C(R_0)}\widehat R(\lambda \pm i0)\chi_{C(R_0)}$.

%%%%%%%%%%%%%%%%%% Lemma 7.5 %%%%%%%%%%%%%%%%%%%%%%
\begin{lemma}
For $\lambda \in I_d\setminus\sigma(\widehat H_{int})$,
$M_{C(R_0)}^{(\pm )} (\lambda )B_{C(R_0)} ^{(\pm )} (\lambda )$ is the identity operator
on $\ell^2(C(R_0))$.
\label{inverse} 
\end{lemma}

%%%%%%%%%%%%%%%%%%%%%%%

\subsection{S-matrix and interior D-N map}

%%%%%%%%%%%%%%%%% Theorem 7.6 %%%%%%%%%%%%%%%%%%%%%

\begin{theorem}
For $\lambda \in I_d \setminus \sigma (\widehat{H}_{int} )$, we have
\begin{equation}
\mathcal{A}_{ext} (\lambda )- \mathcal{A} (\lambda ) =\widehat{ \Gamma }^{(+)} (\lambda ) M^{(+)}_{C(R_0)} (\lambda ) \widehat{\Gamma }^{(-)} (\lambda )^{\ast} . \label{operatoreq} 
\end{equation} 
As a consequence, $ \mathcal{S}(\lambda )$ and $\Lambda _{\widehat{V} } (\lambda ) $  determine each other.
\end{theorem}
Proof. Let us show (\ref{operatoreq}).
For any $\widetilde{\phi } \in L^2 (S^{d-1 })$, let 
\begin{equation}
\begin{split}
 \widehat{u} & = \widehat{ \mathcal{G} }^{(-)} (\lambda )^{\ast } \widetilde{\phi } - \widehat{\mathcal{G}} ^{(-)}_{ext} (\lambda )^{\ast } \widetilde{\phi } \\
 & = \widehat{R} (\lambda +i0 ) \big( \chi_{C(R_0)} B_{C(R_0)}^{(+)}(\lambda ) \chi_{C(R_0)} -\widehat{V} \big) \widehat{\mathcal{G}}_0 (\lambda )^{\ast } \widetilde{\phi }. 
 \end{split}\label{outgoing1} 
 \end{equation}
In view of Lemma \ref{S7exteriorfourier}, $\widehat{u} $ is the outgoing solution of the equation \begin{equation*}
( \widehat{H}_{ext} - \lambda )\widehat{u} =0 \quad \text{in} \quad \stackrel{\circ}\Omega_{ext} , \quad \widehat{u}|_{ \partial \Omega_{ext} } =\widehat{\mathcal{G}}^{(-)} (\lambda ) ^{\ast } \widetilde{\phi } .  
\end{equation*}
By (\ref{S6upmext(n)=RBf}), we can rewrite $\widehat{u} $ as 
\begin{equation}
\widehat{u} = \widehat{R} (\lambda +i0 ) \chi_{C(R_0)} B_{C(R_0)}^{(+)}(\lambda )\chi_{C(R_0)} \widehat{\mathcal{G}}^{(-)} (\lambda )^{\ast } \widetilde{\phi } . \label{outgoing2} 
\end{equation}
By (\ref{outgoing1}), we have as $|k| \rightarrow \infty 
$ \begin{gather*} \begin{split}
\widehat{u}(k)  
&= C_+ |k|^{-(d-1 )/2} e^{ ik \cdot x_{\infty } (\lambda , \omega _k )} a_+ (\lambda , \omega _k ) \\
& \ \ \ \ \ \ \ \ \ \ \ \ \ \times \big( \widehat{\mathcal{G}}^{(+)} (\lambda ) \big(\chi_{C(R_0)} B_{C(R_0) } ^{(+)} (\lambda )\chi_{C(R_0)} -\widehat{V} \big) \widehat{\mathcal{G}}_0 (\lambda )^{\ast }\widetilde{\phi } \big) (\omega_k ) \\
&\quad +O(|k|^{-(d+1 )/2} ), 
\end{split} \end{gather*}
where $C_+ = e^{(3-d) \pi i /4 } \sqrt{2 \pi }$.
On the other hand, by (\ref{outgoing2}), we have as $|k| \rightarrow \infty $ 
\begin{gather*} 
\begin{split}
\widehat{u}(k) 
&= C_+ |k|^{-(d-1 )/2} e^{ ik \cdot x_{\infty} (\lambda , \omega _k )} a_+ (\lambda , \omega _k ) \\
&\ \ \ \ \ \ \ \ \ \ \ \ \ \times \big(\widehat{\mathcal{G}}^{(+)} (\lambda )\chi_{C(R_0)} B_{C(R_0) }^{(+)} (\lambda ) \chi_{C(R_0)} \widehat{\mathcal{G}}^{(-)} (\lambda )^{\ast }\widetilde{\phi } \big) (\omega_k ) \\
&\quad +O(|k|^{-(d+1 )/2}) . 
\end{split} 
\end{gather*}
These two expansions imply 
\begin{equation*}
\begin{split}
 & \widehat{\mathcal{G}}^{(+)} (\lambda ) \big( \chi_{C(R_0)} B_{C(R_0)} ^{(+)} (\lambda )\chi_{C(R_0)} -\widehat{V} \big) \widehat{\mathcal{G}}_0 (\lambda )^{\ast }\\
 & = \widehat{\mathcal{G}}^{(+)} (\lambda ) \chi_{C(R_0)}B_{C(R_0) }^{(+)} (\lambda )\chi_{C(R_0)} \widehat{\mathcal{G}}^{(-)} (\lambda )^{\ast } .
 \end{split} 
 \end{equation*}
The left-hand side is equal to $ \mathcal{A}_{ext} (\lambda )- \mathcal{A} (\lambda )$. On the right-hand side, we insert
$$
 1= M^{(+)}_{C(R_0)}(\lambda)B^{(+)}_{C(R_0)}(\lambda ) :\ell^2 (C(R_0) )\rightarrow \ell^2 (C(R_0))
$$
after $B^{(+)}_{C(R_0)}(\lambda)$
to obtain
\begin{gather*} \begin{split} 
&\widehat{\mathcal{G}} ^{(+)} (\lambda ) \chi_{C(R_0)} B_{C(R_0) }^{(+)} (\lambda ) M^{( +)}_{C(R_0)}(\lambda ) B^{(+ )}_{C(R_0)}(\lambda )  \chi_{C(R_0)} \widehat{\mathcal{G}}^{(-)} (\lambda ) ^{\ast} \\
&=\widehat{\mathcal{G}} ^{(+)} (\lambda ) \chi_{C(R_0)}B_{C(R_0) }^{(+)} (\lambda ) M^{( +)}_{C(R_0)}(\lambda ) \left(B^{(-)}_{C(R_0) }(\lambda )\right)^{\ast } 
\chi_{C(R_0)}\widehat{\mathcal{G}}^{(-)} (\lambda ) ^{\ast} \\
&= \widehat{\Gamma} ^{(+)} (\lambda ) M^{( +)}_{C(R_0)} (\lambda ) \widehat{\Gamma} ^{(-)} (\lambda )^{\ast} . 
\end{split} 
\end{gather*}
We have thus proven (\ref{operatoreq}). 

Given $\Lambda _{\widehat{V} } (\lambda ) $, we can construct $B^{(+)}_{C(R_0)}(\lambda)$ by (\ref{DefineBOmega}), hence $M^{(+)} _{C(R_0)} (\lambda ) $ by Lemma \ref{inverse}.
Since $\widehat{\Gamma }^{(\pm )} (\lambda )$ does not depend on $\widehat{V} $ by Lemma \ref{exteriorasymp},  we can then construct $\mathcal{A}(\lambda )$ by (\ref{operatoreq}).

 By (\ref{operatoreq}), we have
$$
\widehat\Gamma^{(+)}(\lambda)^{\ast}\left(\mathcal A_{ext}(\lambda) - \mathcal A(\lambda)\right)\widehat\Gamma^{(-)}(\lambda) = 
\widehat\Gamma^{(+)}(\lambda)^{\ast}\widehat\Gamma^{(+)}(\lambda)M^{(+)}_{C(R_0)}\widehat\Gamma^{(-)}(\lambda)^{\ast}\widehat\Gamma^{(-)}(\lambda).
$$
 Lemma \ref{neartofarprop} implies $\widehat\Gamma^{(\pm)}(\lambda)^{\ast}\widehat\Gamma^{(\pm)}(\lambda)$ is 1 to 1 on the finite dimensional space $\ell^2(C(R_0))$, hence bijective.  Therefore, one can construct $M^{(+)}_{C(R_0)}(\lambda)$ from $\mathcal{A}(\lambda )$. \qed

%%%%%%%%%%%%%%%%%%%%% Section 8 %%%%%%%%%%%%%%%%%%%

\section{Reconstruction from the D-N map}

In this section, we reconstruction $\widehat V$ from the D-N map $\Lambda_{\widehat V}(\lambda)$.

%%%%%%%%%%%%%%%%%%%%%%%%% subsection 8.1 %%%%%%%%%%%%%%%%%%%%

\subsection{Some properties of Schr\"{o}dinger matrices}
We identify $ -\Delta_{disc} $ and $\Lambda _{\widehat{V}}(\lambda)$  with  matrices as follows. 
Let $ n^{(1)} , \cdots , n^{(\nu)} $ are vertices in $\stackrel{\circ}\Omega_{int} $ and $n^{(\nu +1)} , \cdots , n^{(\nu + \mu )} $ are those in $\partial \Omega_{int} $. We put
$$
\mathcal{N}_0 = \{ n^{(1)} , \cdots , n^{(\nu )} \}, \quad 
\mathcal{N}_1 =\{ n^{(\nu +1)} , \cdots , n^{( \nu+ \mu)} \},
$$
and
$$
\widetilde{\mathrm{deg} }_{\Omega _{int}} (n)=
\left\{
\begin{array}{ll}
^{\#} \{ m\in \Omega_{int} \ ; \ m \sim n \}=2d , & n\in \, \stackrel{\circ}\Omega_{int} , \\
^{\#} \{ m \in  \, \stackrel{\circ}\Omega_{int} \ ; \ m\sim n \} =1, & n\in \partial \Omega_{int} .
\end{array}
\right.
$$
In view of the Laplacian on graphs, we construct a $ (\nu+\mu ) \times (\nu+\mu ) $ matrix ${\bf H}_0 = ( h_{ij}^0 )$ as follows (For the definition, see also \cite{Du84}). 
\begin{gather*} 
{\bf H}_0 =\frac{1}{4} ({\bf D} - {\bf A} ) , \\
{\bf D} = (d_{ij} ), \quad d_{ij} = \left\{ 
\begin{array}{cl} \widetilde{\mathrm{deg} }_{\Omega_{int}} (n^{(i)} ) & (i=j) \\ 0 & ( i\not= j ) 
\end{array} \right. , \\
{\bf A} = (a_{ij} ), \quad
a_{ij} = \left\{ 
\begin{array}{ll} 1, & {\rm if}\quad n^{(i)} \sim n^{(j)} \ \text{for}   \ n^{(i)} \in \, \stackrel{\circ}\Omega_{int} \ \text{or} \ n^{(j)} \in \, \stackrel{\circ}\Omega_{int}, \\ 
0, &  {\rm if}\quad n^{(i)} \not\sim n^{(j)}  , \ \text{or} \ n^{(i)} , n^{(j)} \in \partial \Omega_{int}. 
\end{array} \right. 
 \end{gather*} 
The potential $\widehat{V} $ is identified the diagonal matrix ${\bf V} = (v_{ij } )$ with
 \begin{equation*}
v_{ij }= \left\{ 
\begin{array}{cl} \widehat{V} ( n^{(i)} ) & ( i=j , \ i \leq \nu ) \\ 0 & ( i \not= j \ \text{or} \ i \geq \nu +1 ) 
\end{array} \right. .
 \end{equation*}
Then $\widehat{H}=\widehat{H}_0 + \widehat{V} $ corresponds to the symmetric matrix ${\bf H}={\bf H}_0 +{\bf V} $.
Moreover,  identifying $\widehat u$ with a vector
$( \widehat u(\mathcal N_0) , \widehat u(\mathcal N_1) ) \in {\bf C}^{\nu+\mu}$, 
the equation
\begin{equation}
( -\Delta _{disc} + \widehat{V} ) \widehat{u} =0 \quad \text{in} \quad \stackrel{\circ}\Omega _{int} , 
\label{S8Notmatrixeq}
\end{equation}
is rewritten as 
\begin{equation}
{\bf H}(\mathcal N_0;\mathcal N_1)\widehat u(\mathcal N_1) + 
{\bf H}(\mathcal N_0;\mathcal N_0)\widehat u(\mathcal N_0) = 0,
\label{dirichletprob7int}
\end{equation}
where by ${\bf H}(\mathcal N_i;\mathcal N_j)$ we mean a matrix of size $\, ^{\#} {\mathcal N_i}\times \, ^{\# } {\mathcal N_j}$. The D-N map $\Lambda_{\widehat V}\widehat u =: \widehat g$ is rewritten as
\begin{equation}
{\bf H}(\mathcal N_1;\mathcal N_1)\widehat u(\mathcal N_1) + 
{\bf H}(\mathcal N_1;\mathcal N_0)\widehat u(\mathcal N_0) = \widehat g(\mathcal N_1).
\label{DiscreteDN}
\end{equation}
Taking into account of the Dirichlet data
\begin{equation}
\widehat{u} |_{ \partial \Omega_{int} } =\widehat{f}, 
\label{S8Boudarydataoriginal}
\end{equation}
the above two equations are rewritten as
\begin{equation}
\left( 
\begin{array}{cc} 
{\bf H} ( \mathcal{N}_0 ; \mathcal{N}_0 ) & {\bf H} ( \mathcal{N}_0 ; \mathcal{N}_1 ) \\ {\bf H} (\mathcal{N}_1 ; \mathcal{N}_0 ) & {\bf H} (\mathcal{N}_1 ; \mathcal{N}_1 ) 
\end{array} \right) \left( 
\begin{array}{c} \widehat{u} ( \mathcal{N}_ 0 ) \\ \widehat{f} (\mathcal{N}_1 ) \end{array} 
\right) = \left( 
\begin{array}{c}   \mathbf{0} \\ \widehat{\phi }(\widehat{f} )
\end{array} \right), \quad
 \widehat{\phi }(\widehat{f} ) := \widehat g(\mathcal N_1).
 \label{dirichletprob7total}
\end{equation}

Assume that zero is not a Dirichlet eigenvalue of $-\Delta_{disc} + \widehat{V} $, which means that if $\widehat u(\mathcal N_1) = 0$ in (\ref{dirichletprob7int}), then $\widehat u(\mathcal N_0) = 0$. 
Hence ${\bf H} (\mathcal{N}_0 ; \mathcal{N}_0 )$ is nonsingular. 
Then by using (\ref{dirichletprob7int}), the D-N map corresponds to the $\mu \times \mu$ matrix 
\begin{equation}
{\bf \Lambda _{\widehat{V}}} \widehat{f}(\mathcal N_1):= {\bf H} (\mathcal{N}_1 ; \mathcal{N}_1 ) \widehat{f}(\mathcal N_1) - {\bf H} ( \mathcal{N}_1 ; \mathcal{N}_0 ) {\bf H} (\mathcal{N}_0 ; \mathcal{N}_0 ) ^{-1} {\bf H} (\mathcal{N}_0 ; \mathcal{N}_1 ) \widehat{f}(\mathcal N_1) . 
\label{S7DNdiscrete}
\end{equation}

To simplify the explanation, we translate $\Omega_{int}$ so that 
\begin{equation}
\stackrel{\circ}\Omega_{int} =\{ n \in {\bf Z}^d \ ; \ 1 \leq n_j \leq M , \ j= 1 , \cdots , d \}
\label{Omegainttranslate}
\end{equation}
 for a positive integer $M $. 
We put 
$$
\partial \Omega_j^+ =\{ n \in \partial \Omega_{int} \ ; \ n_j = M+1\}, \quad 
\partial \Omega_j^- =\{ n \in \partial \Omega_{int} \ ; \ n _j =0 \}, \quad j=1, \cdots, d.
$$

%%%%%%%%%%%%%%%%%%%%%% Lemma 8.1 %%%%%%%%%%%%%%%%%%%%

\begin{lemma}
Given a partial Dirichlet data $\widehat{f}$ on $\partial \Omega_{int} \setminus \partial \Omega_1^+$ and a partial Neumann data $ \widehat{g}$ on $\partial \Omega_1 ^-$, there is a unique solution $\widehat{u} $ on $\stackrel{\circ}\Omega_{int} \cup\, \partial \Omega_1^+$ to the equation 
\begin{equation}
\left\{
\begin{split}
& (- \Delta_{disc} + \widehat V)\widehat u = 0 \quad {\rm in}\quad \stackrel{\circ}\Omega_{int}, \\
&\widehat u = \widehat f \quad {\rm on} \quad \partial\Omega_{int}\setminus\partial\Omega_1^+, \\
& \partial_{\nu}^{\Omega_{int}}\widehat u = \widehat{g} \quad {\rm  on} \quad \partial\Omega_1^-.
\end{split}
\right.
\label{S8partialdataequation}
\end{equation}
\label{lemextend}
\end{lemma}
Proof. 
From the boundary values $\widehat{f} ( 0 , n_2 , \cdots , n_d)$ and $\widehat{g}(0 , n_2, \cdots, n_d ) $, we can determine uniquely $\widehat{u} ( 1, n_2 , \cdots , n_d)$ for all $1\leq n_j \leq M$ for $j=2, \cdots , d$: 
\begin{equation*}
\widehat{u} (1 , n_2 ,\cdots , n_d )=- 4 \, \widehat{g} (0,n_2, \cdots , n_d ) + \widehat{f} (0, n_2 ,\cdots , n_d ) .
\end{equation*}
From the equality $((-\Delta_{disc} +\widehat{V} )\widehat{u} )(1, n_2 , \cdots ,n_d )=0$ and the Dirichlet data $ \widehat{f}|_{ \partial \Omega_j^{\pm }}$ for $j=2, \cdots ,d$, we can compute $\widehat{u} (2, n_2, \cdots ,n_d )$ as follows: 
\begin{gather*} 
\begin{split}
\frac{1}{4} \widehat{u} (2, n_2 ,\cdots  , n_d )&= -\frac{1}{4}\sum_{j=2}^d  \sum_{ \alpha =\pm1} \widehat{u} ( 1 , n_2 , \cdots , n_j +\alpha , \cdots , n_d ) -\frac{1}{4} \widehat{f} (0, n_2 ,\cdots , n_d ) \\
&\quad + \frac{d}{2} \widehat{u} ( 1, n_2,\cdots , n_d ) + \widehat{V} ( 1, n_2 ,\cdots , n_d ) \widehat{u} (1, n_2 , \cdots , n_d ), 
\end{split} 
\end{gather*}
for all $ 1 \leq n_j \leq M$, $j=2,\cdots ,d$. 
We repeat this procedure to compute $\widehat{u} (n )$ for all $ n_1 =1 , \cdots,   M+1$. \qed 

\medskip
For subsets $A, B\subset \partial\Omega_{int}$, we denote the associated submatrix of ${\bf{\Lambda}_{\widehat{\bf V}}}$ by ${\bf{\Lambda}}_{\widehat{\bf V}}(B;A).$

%%%%%%%%%%%%%% Corollary 8.2 %%%%%%%%%%%%%

\begin{cor}\label{S8UniqunessCor}
Let $\widehat u$ be the solution of (\ref{S8Notmatrixeq}), (\ref{S8Boudarydataoriginal}). If $\widehat f = 0$ on $\partial\Omega_{int}\setminus\partial\Omega_1^+$, $\Lambda_{\widehat V}\widehat f = 0$ on $\partial\Omega_1^-$, then $\widehat u = 0$ in $\Omega_{int}$.
\end{cor}

%%%%%%%%%%%%%%%%%%%%%%% Corollary 8.3 %%%%%%%%%%%%%%%%

\begin{cor}
The submatrix ${\bf \Lambda _{\widehat{V} }} ( \partial \Omega_1^- ; \partial \Omega_1^+ )$ is nonsingular i.e. ${\bf \Lambda _{\widehat{V}}} (\partial \Omega_1^- ; \partial \Omega_1^+ ) : \partial \Omega_1^+ \rightarrow \partial \Omega_1^- $ is a bijection.
\label{nonsingdnsub} 
\end{cor}
Proof. Suppose $\widehat{f} =0 $ on $\partial \Omega_{int} \setminus \partial \Omega_1^+$ and ${\bf \Lambda _{\widehat{V}}} \widehat{f} =0 $ on $\partial \Omega_1^-$.
By Corollary \ref{S8UniqunessCor}, the solution $\widehat{u} $  of (\ref{S8Notmatrixeq}), (\ref{S8Boudarydataoriginal}) vanishes identically.
Hence $\widehat{f} =0 $ on $\partial \Omega_1^+ $.
This implies that ${\bf \Lambda _{\widehat{V} }} ( \partial \Omega_1^- ; \partial \Omega_1^+ )$ is nonsingular. \qed 

%%%%%%%%%%%%%%%%%%%%% Corollary 8.4 %%%%%%%%%%%%%%%%%

\begin{cor}
Given D-N map ${\bf \Lambda} _{\widehat{\bf V} }$, partial Dirichlet data $\widehat{ f}_2$ on $\partial \Omega_{int} \setminus \partial \Omega_1^+$ and partial Neumann data $\widehat g$ on $\partial\Omega_1^-$, there exists a unique $\widehat f$ on $\partial\Omega_{int}$ such that $\widehat f = \widehat f_2$ on 
$\partial \Omega_{int} \setminus \partial \Omega_1^+$ and
${\bf \Lambda _{\widehat{V} }} \widehat{f} |_{\partial \Omega _1^- } = \widehat g$ on $\partial\Omega_1^-$.
\label{extension} 
\end{cor}
Proof. We seek $\widehat f$ such that
\begin{equation*}
{\bf \Lambda_{ \widehat{V}}} \widehat{f} |_{\partial \Omega_1^- } = {\bf \Lambda _{\widehat{V} }} (\partial \Omega_1^- ; \partial \Omega_1^+ ) \widehat{f} _1 + {\bf \Lambda _{\widehat{V}}} ( \partial \Omega_1 ^- ; \partial \Omega_{int} \setminus \partial \Omega_1^+ ) \widehat{f}_2 = \widehat g,
 \end{equation*}
where $ \widehat{f}_1 = \widehat{f} |_{ \partial \Omega_1^+ } $.
By Corollary \ref{nonsingdnsub}, we take 
$$
\widehat{f}_1 = ({\bf \Lambda _{\widehat{V}}} (\partial \Omega_1^- ; \partial \Omega_1^+ ))^{-1 } \left(\widehat{g}  - {\bf\Lambda _{\widehat{V}} }(\partial\Omega_1^-;\partial \Omega_{int} \setminus \partial \Omega_1^+) \widehat{f}_2 \right). \qed
$$

%%%%%%%%%%%%%%%%%%% subsection 8.2 %%%%%%%%%%%%%%%%%

\subsection{Reconstruction procedure from ${\bf \Lambda _{\widehat{V}}}$}
We can now reconstruct $\widehat{V} $ from ${\bf \Lambda _{\widehat{V} }}$. 
When $d=2$, the procedure has been already given in \cite{Cu1}, \cite{Cu2}, \cite{Ob}. For $d\geq 3$, we generalize this method as follows.

\medskip

\begin{figure}[h]
\centering
\includegraphics[width=6cm, bb=0 0 231 217]{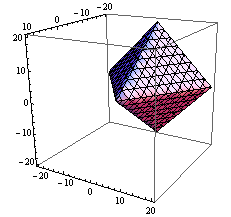}
\caption{The shape of $C_1 (0)$ in the case $d=3$.}
\label{fig:math8}
\end{figure}
We introduce the cone with vertex $n \in \Omega_{int} $ by
\begin{equation}
C_1 (n)=\Big\{ m\in \Omega_{int} \ ; \ \sum_{k\not= 1} | m_k - n_k | \leq - (m_1 - n_1 ) \Big\} .
\label{S8cone}
\end{equation}
If $\widehat{u}$ satisfies the equation (\ref{S8partialdataequation}), we have 
\begin{equation}
\widehat{u} (n ) = \sum_{ m \in C_1 (n) \setminus \{n \}} c_m \widehat{u} (m)
\label{S8continuation}
\end{equation}
for some constants $c_m$.
In particular, if $\widehat{u} (m)=0$ for all $ m\in C_1 (n) \setminus \{n \} $, we see that $\widehat{u} (n)=0$ from (\ref{S8continuation}) (See also Figure \ref{fig:math8}).

\begin{figure}[t]
\centering
\includegraphics[width=9cm, bb=0 0 404 362]{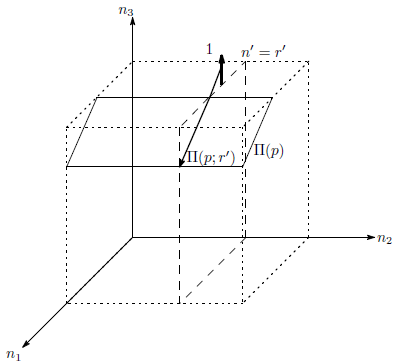}
\caption{Situation of Lemma \ref{S8continuation2}.}
\label{S8figuresituation}
\end{figure}

Let $\Pi(p)$ be the rectangular domain defined by
\begin{equation}
\Pi(p)= \big\{ (n_1,\cdots,n_d) \in \Omega_{int} \, ; \, n_1 + n_d = p, \ 1\leq n_i\leq M \ (2\leq i \leq d-1) \big\},
\label{SetPi(p)}
\end{equation}
where $M$ is from (\ref{Omegainttranslate}),
and for $r' = (r_2,\cdots,r_{d-1}) \in [1,M]^{d-2}$, we consider its section
\begin{equation}
\Pi(p;r') = \big\{ (n_1,n', n_d) \in \Pi(p) \, ; \, n' = r' \big\}.
\label{SectionofPi(p)}
\end{equation}
For $d=3 $, see Figure \ref{S8figuresituation}.

%%%%% Lemma8.5 %%%%%%%
\begin{lemma}
Assume $M+1 < p \leq 2M$, and take a point 
$(p-M-1,r',M+1) \in \Pi(p;r')$. 
Let $\widehat u$ be the solution of (\ref{S8partialdataequation}) with 
 Dirichlet boundary data $\widehat{f} $ such that 
\begin{equation*}
\left\{
\begin{split}
& \widehat{f} (p-M-1,r',M+1) =1 ,\\
& \widehat{f} =0 \quad \text{on} \quad \partial \Omega_{int} \setminus \big( \partial \Omega_1^+ \cup \{(p-M-1,r',M+1)\}\big) ,
\end{split}
\right.
\end{equation*} 
and Neumann data $\widehat g = 0$ on $\partial \Omega_1^-$.
Then we have 
\begin{gather} 
\left\{
\begin{split}
& \widehat{u}(n)= 0   \quad \text{if} \quad n_1 + n_d < p, \\
& \widehat{u}(n)= 0   \quad \text{if} \quad  n_1 + n_d = p, \ n' \neq r', 
\end{split}
\right. 
\label{S8zero}\\
\widehat{u}(p-M-1+i,r',M+1-i)= 
 (-1)^i \quad \text{for}\quad p-M-1+i \leq M+1.
 \label{S8nonzero}  
\end{gather}
If $p=M+1$, taking the Dirichlet data $\widehat{f} $ such that 
\begin{equation*}
\left\{
\begin{split}
&\widehat{f} (0, r', M)=1 ,  \\
&\widehat{f}=0\quad \text{on} \quad \partial \Omega_{int} \setminus \big( \partial \Omega_1^+ \cup \{(0, r' , M) \}\big) ,
\end{split}
\right.
\end{equation*}
we have the same assertion.
\label{S8continuation2}
\end{lemma}

\medskip 

Proof. We put $m= (p-M-1 , r', M+1)$.
First we show that $m\not\in C_1 (n)$, if $n_1 + n_d <p $. In fact,
\begin{equation}
-(m_1-n_1) = n_1 - (p-M-1) < p -n_d - (p-M-1) = M+1-n_d,
\nonumber
\end{equation}
and on the other hand,
\begin{equation}
\sum_{k\neq1}|m_k - n_k| \geq |m_{d} - n_d| = M+1-n_d.
\nonumber
\end{equation}
Then, in view of the condition for $\widehat{f} $, the Neumann data $\partial _{\nu} \widehat{u} |_{\partial \Omega_1^-} =0 $ and (\ref{S8continuation}), we have $\widehat{u} (n)=0 $ if $n_1 + n_d < p$.

Assume that $n_1 + n_d = p$ and $n' \neq r'$.
\begin{figure}[b]
\centering
\includegraphics[width=8cm, bb=0 0 266 150]{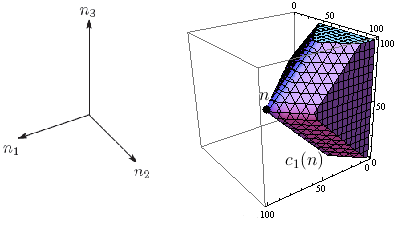}
\caption{Extension of the solution for the case (\ref{S8zero}).}
\label{fig:math10}
\end{figure}
(See Figure \ref{fig:math10}.)
Then 
\begin{equation*}
-( m_1 -n_1 )= M+1-n_d .
\end{equation*}
On the other hand, since $n'\neq r'$, we see that 
\begin{equation*}
\sum_{k\not= 1}|m_k -n_k| > |m_d-n_d| = M+1-n_d .
\end{equation*}
They imply $m \not\in C_1 (n)$, hence $\widehat{u} (n)=0 $ as above.

\begin{figure}[t]
\centering
\includegraphics[width=9cm, bb=0 0 474 278]{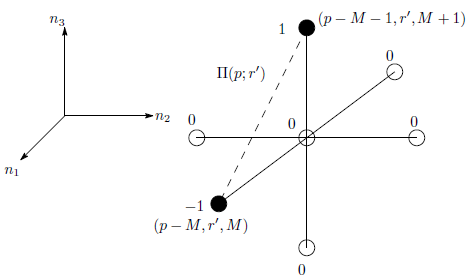}
\caption{Extension of the solution for the case (\ref{S8nonzero}).}
\label{fig:math9}
\end{figure}

Let us prove (\ref{S8nonzero}).
Using the equation 
$$
((-\Delta_{disc} +\widehat{V})\widehat{u})(p-M-1 , r', M) =0,
$$
 and the fact that 
$$
\widehat{u} =0 \quad {\rm for} \quad n_1 + n_d <p, \quad
\widehat u(p-M-1,r',M+1) = 1,
$$
we have $\widehat{u} (p-M , r', M)=-1$. Here we do not use the value of the potential  $\widehat{V}(p-M , r' , M)$. (See Figure \ref{fig:math9}.) 
Repeating this procedure, 
we see $\widehat{u} (p-M-1 +i , r' , M+1-i ) =(-1)^i$ inductively. \qed

\medskip

Now we show the reconstruction procedure.

\medskip

\noindent
{\bf 1st step.} We construct the boundary data $\widehat{f} $ such that 
\begin{equation*}
\left\{
\begin{split}
& \widehat{f} (M-1 , r' , M+1 )=1 , \\
& \widehat{f} =0 \quad \text{on} \quad \partial \Omega_{int} \setminus \big( \partial \Omega_1^+ \cup \{(M-1 , r' , M+1 )\} \big) , \\
& {\bf \Lambda _{\widehat{V}}} \widehat{f}=0 \quad \text{on} \quad \partial \Omega _1^- , 
\end{split}
\right.
\end{equation*}
by Corollary \ref{extension}.
Then the solution $\widehat{u} $ of (\ref{S8Notmatrixeq}) and (\ref{S8Boudarydataoriginal}) satisfies the assumption of Lemma \ref{S8continuation2}.
By virtue of Lemma \ref{S8partialdataequation}, we have 
\begin{equation*}
\widehat{u}(n)= \left\{
\begin{array}{cl}
-1 & (n = (M,r' ,M) ),\\
0 & ( \text{other} \ n\in \, \stackrel{\circ}\Omega_{int} ) .
\end{array}
\right.
\end{equation*}
Then, using the equality 
$$
((-\Delta_{disc}+\widehat{V})\widehat{u}) (M,r' ,M) =0
$$
and the boundary value $\widehat{f} (M+1 , r' ,M) $, we can compute the value $\widehat{V}(M,r' ,M) $.
Applying this procedure for all $r'$,  we recover $\widehat{V} $ on all 
vertices $(n_1,r',n_d)$ such that $n_1 +n_d =2M $.

\medskip

\noindent
{\bf 2nd step.} Assume that we have recovered $\widehat{V}$ on vertices such that $n_1 +n_d >p$ for $ M+1 <p \leq 2M$.
We construct the boundary data $\widehat{f} $ such that 
\begin{equation*}
\left\{
\begin{split}
& \widehat{f} (p-M-1 , r' , M+1 )=1 , \\
& \widehat{f} =0 \quad \text{on} \quad \partial \Omega_{int} \setminus \big( \partial \Omega_1^+ \cup \{(p-M-1 , r' , M+1 )\} \big) , \\
& {\bf \Lambda _{\widehat{V}}} \widehat{f}=0 \quad \text{on} \quad \partial \Omega _1^- . 
\end{split}
\right.
\end{equation*}
By the same argument as in Step 1, the solution $\widehat{u}$ of  (\ref{S8Notmatrixeq}) and (\ref{S8Boudarydataoriginal}) satisfies 
(\ref{S8zero}), (\ref{S8nonzero}).
Since we have already recovered $\widehat{V} $ on $n_1 + n_d >p $, we can compute $\widehat{u}(n) $ on $n_1 +n_d >p$ using the equation $(-\Delta_{disc}+ \widehat{V} )\widehat u= 0$ and the boundary data $\widehat f$.
Hence, using the equality 
$$
((-\Delta_{disc}+\widehat{V})\widehat{u}) (p-M-1+i ,r',M+1-i ) )=0,
$$
and the fact that $\widehat u(p-M-1+i ,r',M+1-i ) = (-1)^{i}$, 
we can compute $\widehat{V}(p-M-1+i ,r' ,M+1-i) $ for every $i$.
Applying this procedure for all $r'$,  we recover $\widehat{V} $ on all 
vertices $(n_1,r',n_d)$ such that $n_1 +n_d =p$ with $M+1<p\leq 2M$.

\medskip

\noindent
{\bf 3rd step.} For $p=M+1$, we construct the boundary data $\widehat{f} $ such that 
\begin{equation*}
\left\{
\begin{split}
& \widehat{f} (0 , r' , M )=1 , \\
& \widehat{f} =0 \quad \text{on} \quad \partial \Omega_{int} \setminus \big( \partial \Omega_1^+ \cup \{(0 , r' , M )\} \big) , \\
& {\bf \Lambda _{\widehat{V}}} \widehat{f}=0 \quad \text{on} \quad \partial \Omega _1^- . 
\end{split}
\right.
\end{equation*}
By the same argument as in Step 1, the solution $\widehat{u}$ of  (\ref{S8Notmatrixeq}) and (\ref{S8Boudarydataoriginal}) satisfies 
\begin{equation*}
\widehat{u}(n)= \left\{
\begin{array}{cl}
(-1)^{i-1} & (n = (i , r' ,M+1 -i ) ),\\
0 & (n_1+n_d < p \ \  {\rm or}\ \ n_1+n_d = p,\  n'\neq r' ) .
\end{array}
\right.
\end{equation*}
Then we can compute $\widehat{V} (i ,r',M+1 -i ) $ for every $i$ as above.

\medskip

\noindent
{\bf 4th step.} In the case $n_1 +n_d <M+1 $, we have only to rotate the whole domain.

\medskip
We have thus completed the proof of Theorem \ref{Main theorem}.

%%%%%%%%%%%%%%%%%%%%%%%%%%%%%%%%%%%%

\end{document}